\documentclass{amsart}

\newtheorem{theorem}{Theorem}[section]

\newtheorem{definition}[theorem]{Definition}
\newtheorem{example}[theorem]{Example}

\theoremstyle{remark}
\newtheorem{remark}[theorem]{Remark}

\numberwithin{equation}{section}

\usepackage{amsmath}
\usepackage{amsfonts}
\usepackage{rotating}
\usepackage{bbm}
\usepackage[all,cmtip]{xy}

\usepackage{amssymb}
\usepackage[latin1]{inputenc}
\usepackage{graphicx}
\usepackage{dsfont}
\usepackage{color}
\usepackage{hyperref}
\usepackage{amssymb,amsthm,amsfonts,amsmath,amscd,amssymb,epsfig}
\usepackage[all]{xy}

\usepackage{enumerate,array,hyperref, graphicx,color,amsthm,
bbm,mathrsfs,a4wide}
\usepackage{amssymb}   
\usepackage{graphicx}
\usepackage{caption}
\usepackage{setspace}
\usepackage{tikz}

\renewcommand{\epsilon}{\varepsilon}

\renewcommand{\S}{S}

\begin{document}
\title{The ECH capacities for the rotating Kepler problem}
\author{Amin Mohebbi}
\keywords{symplectic geometry, special concave toric domain, the new tree, ECH capacities, embedding problem, the rotating Kepler problem, regularization problem. }
\let\thefootnote\relax\footnote{Partially funded by FCT/Portugal through UID/MAT/04459/2020 and project PTDC/MAT-PUR/29447/2017.} 

\maketitle
\section{Introduction}
In this paper, I am going to use the Special Concave Troic Domain of the Rotating Kepler Problem \cite{key-1} and compute the ECH capacities of the RKP for all energy less than or equal to the critical energy value, i.e. for all energy $c\leq -\dfrac{3}{2}$. To compute these ECH capacities, first, we need to obtain the weights of the SCTD for the RKP, then find an order for the weights from the biggest to the lowest one.  
  
From \cite{key-1}, we know that the SCTD is a Concave Toric Domain after applying the Ligon-Schaaf symplectomorphism \cite{key-5}  and Levi-Civita  regularization \cite{key-14} on the CTD.  Thus we can say the SCTD is a special case of the CTD which is rotated by 45 degree in the clockwise direction. The SCTD gives us a family of CTD's such that the energy parameter is less than or equal to the critical energy value. 
  
For Theorem A, we assume the energy $c \leq - \dfrac{3}{2}$ and try to obtain the weights  to compute the ECH capacities of the RKP and explain all statements of the weight extension and also the method of getting the weights via the new tree \cite{key-1}. For Theorem B, using the concepts and notations in Theorem A, we discuss the ECH capacities for the critical energy, i.e. $c=-\dfrac{3}{2}$  and we will prove that the first weight is the biggest weight for all energy in the SCTD. Finally, we can see a numerical example of the ECH capacities.

\section{Introduction to ECH capacities}
ECH capacities give obstructions to symplectic embeddings of one symplectic 4-manifold with boundary into another one. 
Here we are going to give the definition and some properties of ECH capacities and a new interpretation of Hutchings algorithm \cite{key-9} with help of   the Stern-Brocot tree  \cite{key-1} and \cite{key-10}. This  is useful later for rotated concave toric domain in rotating Kepler problem.

\begin{definition}
Suppose $(X,\omega)$  is a compact symplectic 4-manifold. This manifold  can have boundary and  corners.
 ECH capacities are defined for the manifold $(X,\omega)$ as a sequence of real numbers 
\begin{align}
0= c_0(X,\omega) \leq c_1 (X, \omega) \leq c_2 (X, \omega) \leq \cdots \leq \infty
\end{align}
which have useful properties  as follows. 
\begin{itemize}
\item[(Monotonicity)]
 If there exists a symplectic embedding 
\begin{align}
(X,\omega) \longrightarrow (X',\omega ').
\end{align}
 Then for all $k$, we have inequality 
\begin{align}
c_k (X,\omega) \leq c_k (X',\omega ').
\end{align}
\item[(Conformality)]
For $r >0$ it holds true that 
\begin{align}
c_k(X, r\omega) = r c_k(X,\omega).
\end{align}
\item[(Disjoint Union)] 
\begin{align}
c_k(\bigsqcup_{i=1}^n (X_i,\omega_i)) =\max_{k_1+k_2+ \cdots + k_n = k} \sum_{i=1}^n c_{k_i} (X_i,\omega_i).
\end{align}
\item[(Ellipsoid)] Let $a,b >0$ and define the ellipsoid by 
\begin{align}
E(a,b):= \{ (z_1,z_2) \in \mathbb{C}^2 | \frac{\pi |z_1|^2}{a} + \frac{\pi |z_2|^2}{b} \leq 1 \}.
\end{align}
\end{itemize}
We can write $c_k(E(a,b)) = N(a,b)$, where $N(a,b)$ denotes the sequence of all nonnegative integer linear combinations of $a$ and $b$ arranged in nondecreasing order and index $k$ starting from zero. 
\end{definition}

Consider the standard symplectic form on $\mathbb{C}^2=\mathbb{R}^4$ and  let $a=b$, we can abbreviate
\begin{align}
E(a,b)=E(a,a) =: B(a)
\end{align}
that is called a ball with radius $\sqrt{\frac{a}{\pi}}$. 

If we apply the identity $a=b$ on the ellipsoid property, we get a similar property for a ball with radius $\dfrac{a}{\pi}$, i.e.
\begin{align}
c_k(B(a))= ad
\end{align}
where $d$ is the unique nonnegative integer such that
\begin{align}
\dfrac{d^2+d}{2} \leq k \leq \dfrac{d^3+3d}{2}.
\end{align}
McDuff showed \cite{key-11} that there exists a symplectic embedding $int(E(a,b)) \hookrightarrow E(a',b')$ if and only if $N(a,b)_k \leq N(a',b')_k$ for all $k$. Therefore, ECH capacities give a sharp obstruction to symplectic embedding one (open) ellipsoid into other one. 

Define the  polydisk  
\begin{align}
P(a,b) = \{ (z_1,z_2) \in \mathbb{C}^2  \: : \:  \pi|z_1|^2 \leq a, \: \pi |z_1|^2 \leq b \}.
\end{align}
We use  ECH capacities and  give a  sharp obstruction which is symplectically embedding by
$$E(a,b) \hookrightarrow^{s} P(a',b').  $$
In genera case, the inverse of the above embedding is not hold \cite{key-12}, i.e. ECH capacities give not a sharp obstruction to embedding $P(a',b')$ into $E(a,b)$. 

Here we are going to  introduce a new method to compute the ECH capacities of a Hutchings concave toric domain. In this method, we will use  the  Stern-Brocot tree to obtain the slopes of each portion of the CTD. Note that in the following, we will extend this method to compute the ECH capacities of the rotating Kepler problem using the Special Concave Toric Domain \cite{key-1} and \cite{key-2}.
To this purpose, first we recall the definition of the Hutchings concave toric domain. 

Suppose $\Omega$ is a domain in the first quadrant of the plane $\mathbb{R}^2$. The toric domain is defined by 
\begin{align}
X_\Omega := \{ z \in \mathbb{C}^2 \;\; | \;\; \pi(|z_1|^2 , |z_2|^2) \in \Omega \}.
\end{align}
The map
\begin{align}
\nu : &X_\Omega \longrightarrow \pi \\
z &\mapsto \pi (|z_1|^2 , |z_2|^2),
\end{align}
is called the momentum map. 
\begin{definition}
A concave toric domain is a domain $X_\Omega$ where $\Omega$ is the closed region bounded by the horizontal segment from $(0,0)$ to $(a,0)$, the vertical from $(0,0)$ to $(b,0)$ and the graph of a convex function $f:[0,a] \longrightarrow [0,b]$ with $f(0)=b$ $f(a)=0$. The concave toric domain $X_\Omega$ is rational if $f$ is piecewise linear and $f'$ is rational wherever it is defined. 
\end{definition}
\begin{example}
Given  a triangle with vertices $(0,0), \: (a,0)$ and $(0,b)$ on the standard  coordinate space in $\mathbb{R}^2 $. This  concave toric domain is an ellipsoid such as  $E(a,b)$. 
\end{example}

\subsection{Weight Expansions:}
Let $X_\Omega$ be a CTD, the weight expansions of $\Omega$ is a finite (or infinite) unordered list of (possibly) repeated positive real numbers $W(\Omega) = ( a_1,a_2, \cdots ,a_n)$ defined inductively, such that  the weight belong to portions of the CTD. 

Now we describe how we can get this list of positive real numbers via the Stern-Brocot tree. 
\begin{remark}
By the Stern-Brocot tree, we can relate every weight to a node of the Stern-Brocot tree. On the other hand, these nodes correspond to tori $T_{k,l}$ and we can find their  slopes by the formula
\begin{align}
S_{k,l}^{CTD} = -  \dfrac{k}{l}.
\end{align} 
\end{remark}
\textbf{Recall:}
 A node of the Stern-Brocot tree is called even or odd if we write it by a sequence of 0 and 1 such that the  sequence ends with 0 or 1 respectively. 

Denote the portions of a CTD with $\Omega_{i_1,i_2,\cdots ,i_j}$ where $i_1, \cdots ,i_j \in \{0,1\}$. To each portion like $\Omega_{i_1,i_2,\cdots ,i_j}$  we related the node $N_{i_1,i_2, \cdots , i_j}$ in the Stern-Brocot tree. Using these  notations help us to give the computation of  the weight expansion of the CTD as follow. 

The easiest case appears when $\Omega$ be a triangle with vertices $(0,0), \: (a,0)$ and $(0,a)$. The weight of $\Omega$ in this case is equal to $a$, i.e. $W^{CTD}(\Omega) = a$. 

Otherwise, let $a>0$ be the largest real number such that the triangle with vertices $(0,0),\: (a,0)$ and $(0,a)$is contained in $\Omega$. We name this triangle $\Omega_1$. We have the torus $T_{1,1}$ and the slope $S_{1,1}^{CTD} =-1$ corresponded to portion $\Omega_1$. Using the momentum map $\nu$, we can write $T_{1,1} = \nu^{-1}(\nu_{1\Omega_1} , \nu_{2\Omega_1})$. Hence the first weight of the CTD is 
\begin{align}
W^{CTD}(\Omega_1) =a.
\end{align}
We draw the line $x+y=a$ and obtain the tangent point of it and the graph of $f$. Denote the tangent point with $(\nu_{1\Omega_1} , \nu_{2\Omega_1})$ and call is the critical point of $\Omega_1$. 

 \begin{center}
\includegraphics[scale=0.25]{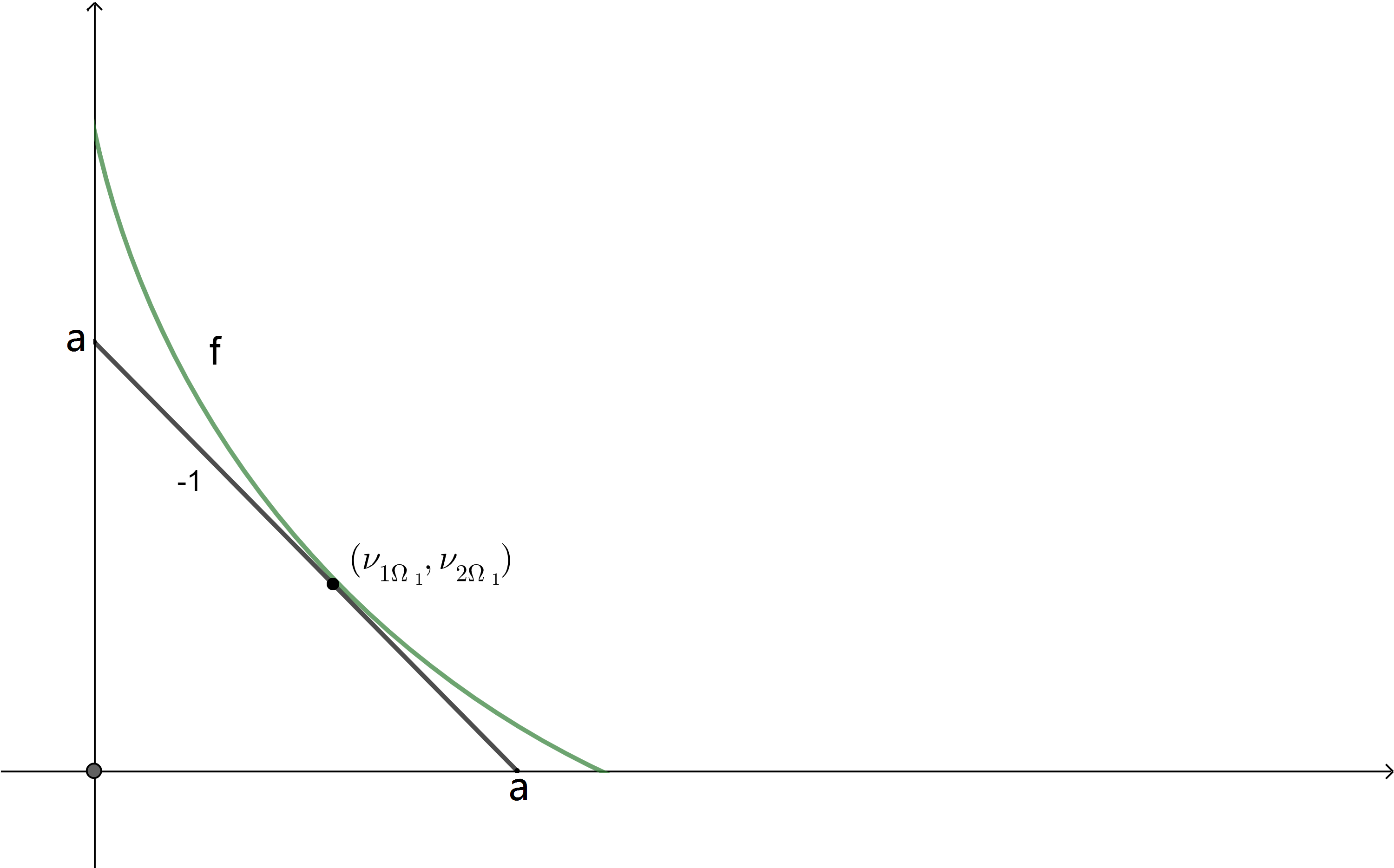}
  \captionof{figure}{ The portion for the first weight $W^{CTD} (\Omega_1)$}
  \label{picCTD5}
\end{center}

For the second weight, we consider the portion $\Omega_{11}$ of $\Omega$ such that this portion lives between the line $x+y=a$ and the graph of the function $f$, above the point $(\nu_{1\Omega_1},\nu_{2\Omega_1})$. Because of the index, we named the portion $\Omega_{11}$ odd. The nodes corresponds to this portion in the Stern-Bocot tree is $\frac{k}{l} = \frac{2}{1}$ and its slope satisfies  
\begin{align}
S_{2,1}^{CTD} = -\dfrac{2}{1} = -2. 
\end{align}
The critical point of the portion $\Omega_{11}$ is $(\nu_{1\Omega_{11}} , \nu_{2\Omega_{11}})$ which is the tangent point of the slope $S_{2,1}^{CTD} = -2$ in the graph of the function $f$. On the other hand, we have 
\begin{align}
T_{2,1} = \nu^{-1} (\nu_{1\Omega_{11}} , \nu_{2\Omega_{11}}).
\end{align}
Denote the intersection point of the slope $S_{2,1}^{CTD}=-2$ and  $y$ axis with $(x_2,y_2)$. The second weight of the CTD is 
\begin{align}
W^{CTD}(\Omega_{11}) = y_2 - W^{CTD}(\Omega_1).
\end{align}
  \begin{center}
\includegraphics[scale=0.22]{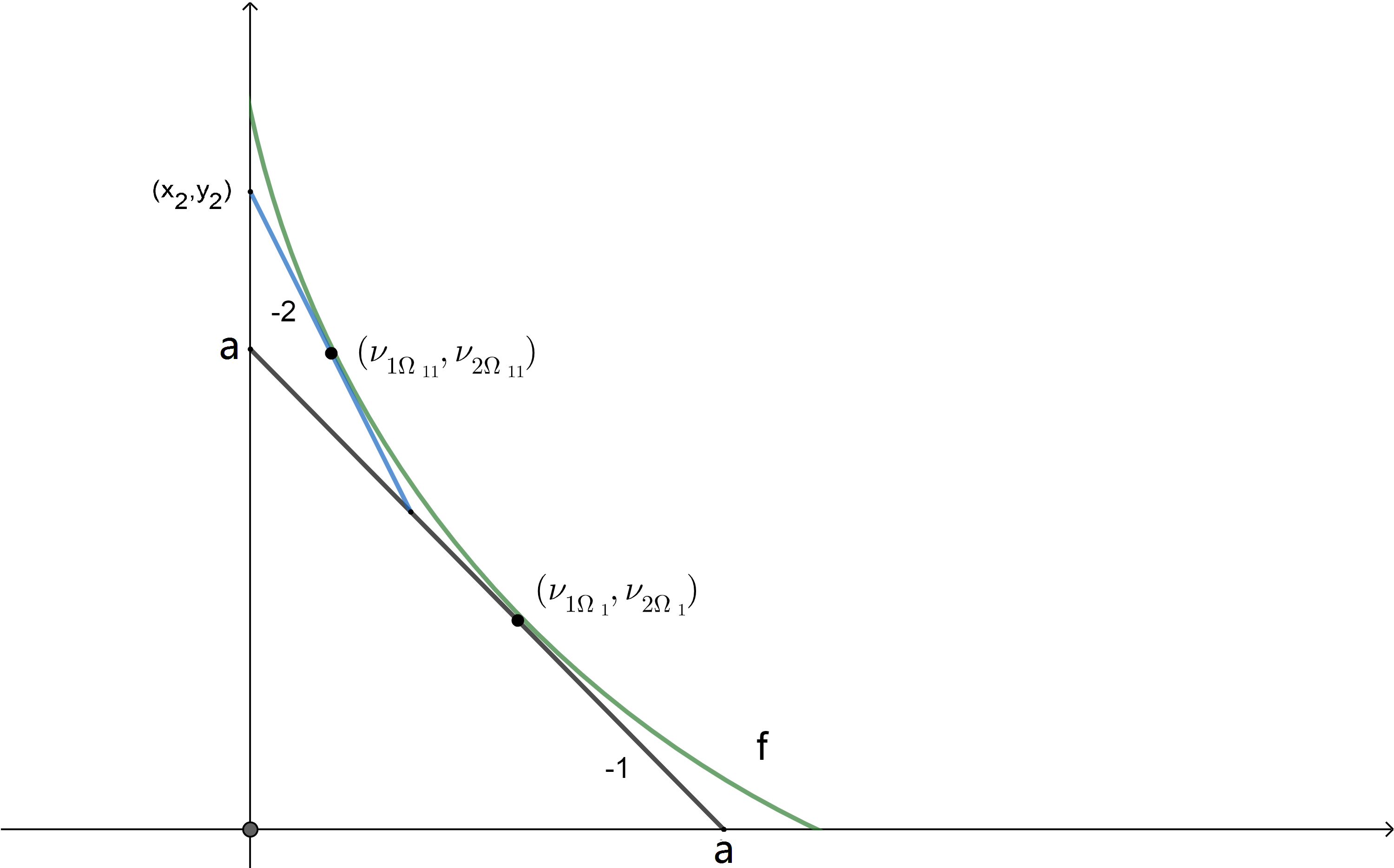}
  \captionof{figure}{ The portion for the  weight $W^{CTD} (\Omega_{11})$}
  \label{picCTD7}
\end{center}

Note that if we use the method in \cite{key-9}, we can convert the portion $\Omega_{11}$ to the standard shape, namely the same shape as $\Omega_2$ by multiplication to the matrix
$
\left[ \begin{array}{ccc}
1 & 0 \\
1 & 1 
\end{array} \right]
 \in SL_2(\mathbb{Z})$ and a transformation. 

Now we assume the portion below the critical point $(\nu_{1\Omega_1} , \nu_{2\Omega_1})$ and denote it by $\Omega_{10}$. According the index, we named $\Omega_{10}$ is an even portion of the CTD $\Omega$. This portion related to the node $N_{1,0}$ in the Stern-Brocot tree which is $\frac{k}{l}=\frac{1}{2}$ and also to the torus 
\begin{align}
T_{1,2} = \nu^{-1}(\nu_{1\Omega_{10}} , \nu_{1\Omega_{10}}),
\end{align}
where $ (\nu_{1\Omega{10}} , \nu_{2\Omega_{10}})$ is the critical point of the portion $\Omega_{10}$ and the slope corresponded to this portion via the Stern-Brocot tree is 
\begin{align}
S_{1,2}^{CTD}=-\dfrac{1}{2}.
\end{align}
If we define the intersection point of the slope $S_{1,2}^{CDT}$ and the x-axis with $(x_3,y_3)$, then we can get the third weight of the CTD as 
\begin{align}
W^{CTD}(\Omega_{10}) = x_3 - W^{CTD}(\Omega_1),
\end{align} 
   \begin{center}
\includegraphics[scale=0.22]{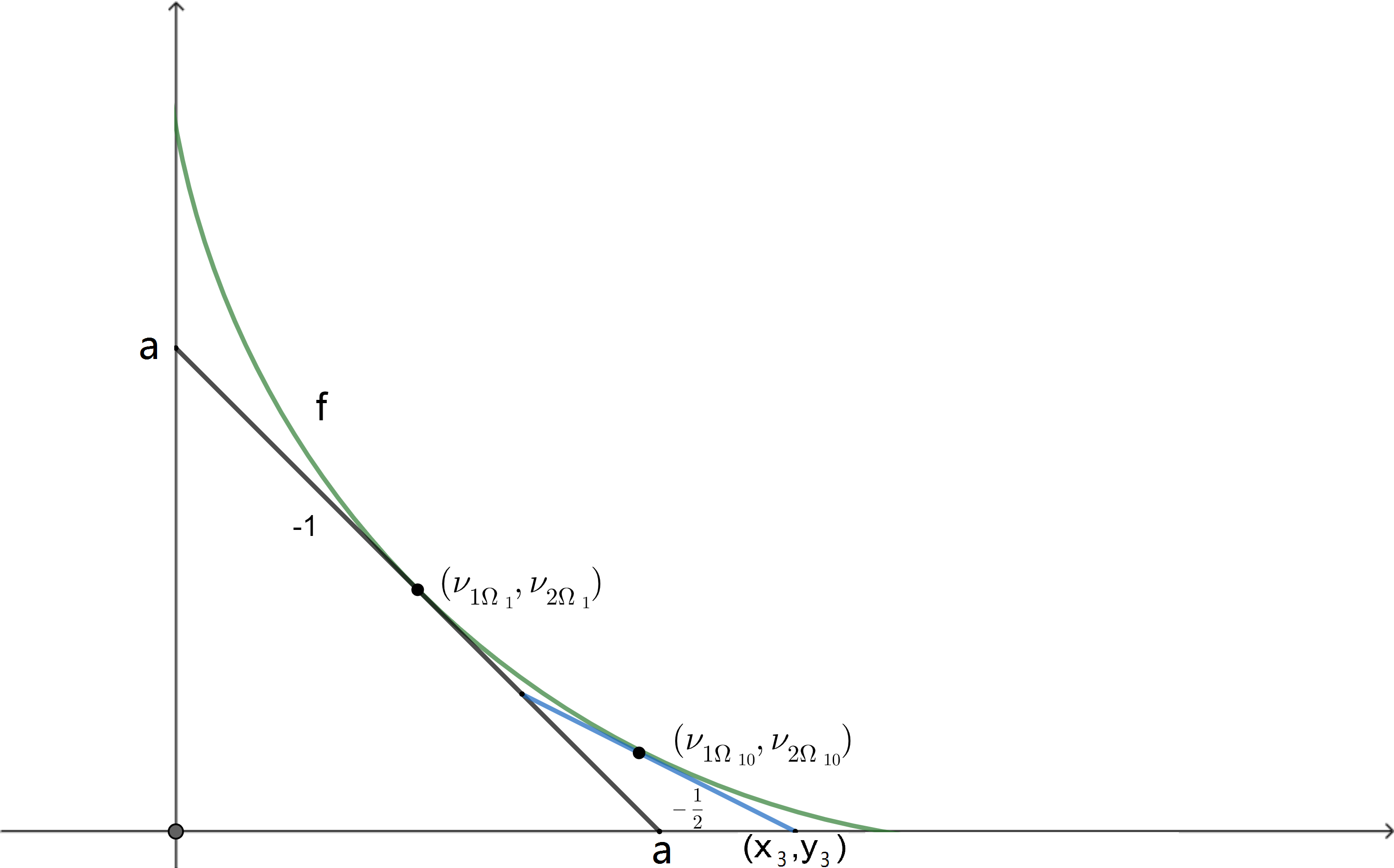}
  \captionof{figure}{ The portion for the  weight $W^{CTD} (\Omega_{10})$}
  \label{picCTD6}
\end{center}

As the portion $\Omega_{11}$, we can convert the portion $\Omega_{10}$ to the standard shape with multiplying to 
 $
\left[ \begin{array}{ccc}
1 & 1 \\
0 & 1 
\end{array} \right]
 \in SL_2(\mathbb{Z})$
 and a transformation. Then we can get a right angle shape for the portion $\Omega_{10}$.

Depending on the even or odd portion of the CTD, we can follow  the above method and compute the higher weights of the CTD belongs to the portion $\Omega_{100}, \Omega_{101}, \Omega_{110}, \Omega_{1000}, \Omega_{1001},   \cdots$. 

After  finding the weights expansion $(a_1,\cdots , a_n)$ for the CTD, we need to give the weights  an order from highest  to lowest weight. To this deal, assume a finite sequence binary number $V$ and define 
\begin{itemize}
\item[(i)]
\begin{align*}
a_1 &:=\max\{W(\Omega_v) \: \: : \: \: v \in V \} \\
n_1&:= \# \{ v \in V \;\; : \;\; W(\Omega_v) = a_1 \}\\
a_i & = a_1 \;\;\; \text{for } 1 \leq i \leq n_1, \\
\end{align*}
\item[(ii)] 
\begin{align*}
a_{n_1+1} &:=\max\{W(\Omega_v) \: \: : \: \: v \in V \; \text{and}\;\; W(\Omega_v) < a_1\} \\
n_{2}&:= \# \{ v \in V \;\; : \;\; W(\Omega_v) = a_{n_1+1} \}\\
a_i & = a_{n_1+1} \;\;\; \text{for } n_1+1 \leq i \leq n_1 + n_2, \\
\end{align*}
recursively 
\item[(iii)]
\begin{align*}
a_{n_k+1} &:=\max\{W(\Omega_v) \: \: : \: \: v \in V \; \text{and}\;\; W(\Omega_v) < a_1\} \\
n_{k+1}&:= \# \{ v \in V \;\; : \;\; W(\Omega_v) = a_{n_k+1} \}\\
a_i & = a_{n_k+1} \;\;\; \text{for } \sum_{j=1}^k  n_j+1 \leq i \leq  \sum_{j=1}^{k+1}n_j. \\
\end{align*}
\end{itemize}

Note that $a_1 \geq a_2 \geq a_3 \geq  \cdots$   and consider $W^{CTD}(\Omega_{i_1i_2\cdots i_j}) =0 $ when $\Omega_{i_1i_2\cdots i_j} =\emptyset $.
\begin{remark}
In the case when $\Omega$ is a rational triangle,  the weight expansion is determined  by continued fraction expansion of the slope of the diagonal and in particular $W(\Omega)$ is finite. 
For instance,  when $X_\Omega$ is rational its weight expansion is finite.
\end{remark}
\begin{theorem}
Given $X_\Omega$  a rational concave toric domain and its weight expansion be $(a_1 , \cdots , a_n)$. The ECH capacities of $X_\Omega$ are given by 
\begin{align}
c_k(X_\Omega) = c_k (\coprod_{i=1}^n B(a_i)). 
\end{align}
\end{theorem} 
\begin{proof}
See \cite{key-9}.
\end{proof}

\section{Computing Some ECH capacities for the Rotating Kepler Problem}
Here, using the Special Concave Toric Domain $\mathcal{K}_c^b$  \cite{key-1}, we are going to  obtain a family of SCTD related to the RKP with energy parameter $c$ and  give the computation of some ECH capacities for this family when the energy $c$ is less than or equal to $-\frac{3}{2}$.  To this purpose, first we recall the SCTD $\mathcal{K}_c^b$  for the RKP. 




Let the critical energy value $-\frac{3}{2}$ and consider the equation  
\begin{align}\label{p1}
\mu_2=\dfrac{1}{16\mu_1^2} +\dfrac{c}{2},
\end{align}
when $c=-\dfrac{3}{2}$.

First we need to find the weights of the SCTD $\mathcal{K}_c^b$ than  compute the ECH capacity of corresponding to each weight. 
\subsection{The First Weight $W_1$}
To compute the first weight, we try to find the intersection point of the graph of equation \ref{p1} and the line $\mu_1=\mu_2$.  If we plug in the equality $\mu_1=\mu_2$ into the equation \ref{p1}. Then we  have 
\begin{align} \label{p2}
-16\mu_1^3+8c\mu_1^2+1=0.
\end{align}
Using the trigonometric  method, we can find the roots of the above equation. To see the details of the roots computation see \cite{key-1}.

The roots of the equation \ref{p2} are
\begin{align}
r_i=( (\frac{c}{3} \cos (\frac{1}{3} \arccos (1+\frac{27}{4c^3}))+\frac{2 \pi}{3} T)) +\frac{c}{6}), \;\;\; i=1,2,3,
\end{align}
where $T=0,1,2$.

 Assume the first root appears when $T=1$, so  we have 
\begin{align}
r_1 (c)= (\frac{c}{3} \cos (\frac{1}{3} \arccos (1+\frac{27}{4c^3})+\frac{2 \pi}{3} )) +\frac{c}{6}
\end{align}
The first weight $W_1$ in the SCTD $\mathcal{K}_c^b$ is the diameter of the isosceles right-angled triangle with the length of the sides $r_1$. Hence, we can use the Pythagorean theorem and write the first weight $W_1$ as a function of $r_1$ by
\begin{align}
W_1(r_1(c))=\sqrt{2}r_1(c), 
\end{align}
or equivalently as a function of the energy $c$ by
\begin{align}
W_1(c) = \sqrt{2}((\frac{c}{3} \cos (\frac{1}{3} \arccos (1+\frac{27}{4c^3})+\frac{2 \pi}{3} )) +\frac{c}{6}).
\end{align}
We can get the roots $r_2$ and $r_3$ as functions of the first root $r_1$ \cite{key-1}, as follows
\begin{align}
r_2&= -\dfrac{-1+\sqrt{1-4^3r_1^3}}{32r_1^2}, \\
r_3& =-\dfrac{-1-\sqrt{1-4^3r_1^3}}{32r_1^2}.
\end{align}
The second root $r_2$ is the intersection point of the graph of the equation \ref{p2} and the line $y=-x$ which plays an important rule to compute the higher weights of the SCTD $\mathcal{K}_c^b$.

The  first weight $W_1$ is a smooth function of $c$. We will see that the higher weights on the SCTD $\mathcal{K}_c^b$  are not smooth, but they are continuous in $c$  in the nonsmooth points. To given these functions, we need to find their domains that they determine by the unique critical energy value for each weight. We will use the new tree \cite{key-1} to find the portions for each weight on the SCTD $\mathcal{K}_c^b$. 

$\textbf{Recall:}$ The new tree is as 

 \begin{center}
\begin{tikzpicture}[<-,->=stealth', auto, semithick,
level/.style={sibling distance=60mm/#1}
]
\node (z){$\infty $}
child {node (a) {$\dfrac{3}{1}$}
child {node (b) {$\dfrac{2}{1}$}
child {node {$\dfrac{5}{3}$}
child {node {$\vdots$}} 
} 
child {node {$\dfrac{7}{3}$}
 child {node {$\vdots$}}}
}
child {node (g) {$\dfrac{5}{1}$}
child {node {$\dfrac{4}{1}$}
child {node {$\vdots$}}}
child {node {$\dfrac{7}{1}$}
child {node {$\vdots$}}}
}
}
child {node (j) {$-\dfrac{3}{1}$}
child {node (k) {$-\dfrac{5}{1}$}
child {node {$-\dfrac{7}{1}$}
child {node {$\vdots$}}}
child {node {$-\dfrac{4}{1}$}
child {node {$\vdots$}}}
}
child {node (l) {$-\dfrac{2}{1}$}
child {node {$-\dfrac{7}{3}$}
child {node {$\vdots$}}}
child {node (c){$-\dfrac{5}{3}$}
child {node {$\vdots$}}
}
}
};
child [grow=down] {
edge from parent[draw=none]
};
(z) -- ++(0,1cm).
\end{tikzpicture}
 \end{center}

This tree helps us to find the slopes  and the critical energies of the portions in the SCTD $\mathcal{K}_c^b$.  

\subsection{The Critical Energy Value and the Slopes of Weights}
In \cite{key-1}, we denote the nodes  $\frac{k}{l}$ of the Stern-Brocot tree related to  the slope $\S_{k,l}^{CTD} =-\frac{k}{l}$ in the Hutchings CTD. Denote the portions of the SCTD $\mathcal{K}_C^B$ by $\omega_{i_1,\cdots , i_j}$ uniquely where $i_n \in \{0,1\}$. We will use this notations to find the weights, slopes and the critical energy values for the SCTD $\mathcal{K}_c^b$ and computing the ECH capacities of the RKP. 

Let the portion $\omega_{i_1,\cdots , i_j}$,  we can compute the critical energy  value  and the slope of the portion  $\omega_{i_1,\cdots , i_j}$ respectively by the following functions, \cite{key-1}, \cite{key-3} 
\begin{align}\label{p3}
c_{k,l}^- &= -\dfrac{1}{2}(\dfrac{k}{l})^{\frac{2}{3}} - (\dfrac{l}{k})^{\frac{1}{3}}\\
&= -(\dfrac{1}{2} +\dfrac{l}{k})(\dfrac{k}{l})^{\frac{2}{3}}
\end{align}
and 
\begin{align}\label{p4}
S_{l,k}=-\dfrac{k+l}{-k+l}.
\end{align}
Note that some tori $T_{k,l}$ are assigned to asteroids. For instance, the tori $T_{2,1}$ and $T_{3,1}$ are assigned to  the asteroids Hekuba and Hestia respectively. 

Note that the critical energy $c_{k,l}^-$ is the energy that the torus $T_{k,l}$ appears first.

Consider a portion $\omega_{i_1,\cdots , i_j}$ for $i_1,\cdots , i_j \in V$ on the SCTD $\mathcal{K}_c^b$. This  portion corresponds to the torus $T_{k,l}$ in the SCTD $\mathcal{K}_c^b$ via the following relation 
\begin{align}\label{p6}
T_{k,l} = \mu^{-1} (\mu_{1\omega_{i_1 \cdots i_j}} , \mu_{2\omega_{i_1 \cdots \i_j}})
\end{align}
where $(\mu_{1\omega_{i_1,\cdots , i_j}} ,\mu_{2\omega_{i_1,\cdots , i_j}})$ is the tangent point of the slope $S_{k,l}$ and the graph of the equation \ref{p1}. We called the point  $(\mu_{1\omega_{i_1,\cdots , i_j}} ,\mu_{2\omega_{i_1,\cdots , i_j}})$,  the critical point of the portion $\omega_{i_1,\cdots , i_j}$.

Note that the tori computed for  portions $\Omega_{i_1,\cdots , i_j}$  by the equation  
  \begin{align} \label{p7}
  T_{k,l} = \nu^{-1} (\nu_{1\Omega_{i_1 \cdots i_j}} , \nu_{2\Omega_{i_1 \cdots i_j}})
  \end{align}
on the CTD are correspond to the tori computed for portion $\omega_{i_1,\cdots , i_j}$ in the SCTD $\mathcal{K}_c^b$. 

Take the first derivative of the equation \ref{p1}. That is 
\begin{align}\label{p5}
\dfrac{d \mu_2}{d \mu_1} = -\dfrac{1}{8 \mu_1^3}. 
\end{align}
If we put the above equation equal to the slope $S_{k,l}$, then we can get the first critical value $\mu_{1\omega_{i_1,\cdots , i_j}}$. 
Now we substitute $\mu_{1\omega_{i_1,\cdots , i_j}}$ into the equation \ref{p1} to get the second critical value $\mu_{2\omega_{i_1,\cdots , i_j}}$. Therefore, the critical point of the portion $\omega_{i_1, \cdots , i_j}$ is $(\mu_{1\omega_{i_1,\cdots , i_j}}, \mu_{2\omega_{i_1,\cdots , i_j}})$. 

Some tori have special name that we will use the special name of them to show their  critical points. In other words, we will show the critical point of the torus $T_{2,1}$ - asteroid Hekuba - and  the critical point of the torus $T_{3,1}$, asteroid Hestia, by $(\mu_{1Hek}, \mu_{2Hek})$  and $(\mu_{1Hes}, \mu_{2Hes})$ respectively. 

\section{The Higher Weights}
\subsection{The second weight $W_2$:}
The second weight $W_2$ belongs to the portion $\omega_{11}$ of the SCTD $\mathcal{K}_c^b$. This portion is the biggest triangle in the rest part of $\mathcal{K}_c^b$ which is bounded by the lines $x=r_1$, $y=-x$ and the graph of the equation \ref{p2}. Using the new tree, we  get the slope $S_{2,1}=-3 $ for the portion $\omega_{11}$. In view of the relation \ref{p4}, the slope of $\omega_{11}$ in $\mathcal{K}_c^b$ is 
\begin{align}
S_{k,l} = S_{2,1} =\dfrac{2+1}{-2+1} =-3,
\end{align}
which is the value of the node $V_{11} $ in the new tree. 

From now, we want to use the new tree to find the slopes of a portion. 

To obtain the weight $W_2$,  we compute the critical energy value using  \ref{p3} as follow
\begin{align}
c_{k,l}^- = c_{2,1}^- = - (\dfrac{1}{2} + \dfrac{l}{k})(\dfrac{k}{l})^\frac{2}{3} = - \sqrt[3]{4}. 
\end{align}
The critical energy value $c_{2,1}^- =-\sqrt[3]{4}$ is the energy which the asteroid Hekuba appears first. 

Now use the relation \ref{p5} to compute the critical point $(\mu_{1\omega_{11}} , \mu_{2\omega_{11}})$ for the portion $\omega_{11}$ in $\mathcal{K}_c^b$  when this relation in equal to the slope $S_{2,1}=-3$,
\begin{align}
\dfrac{d \mu_{2\omega_{11}}}{d \mu_{1\omega_{11}}} =-\dfrac{1}{8\mu_{1\omega_{11}^3}}.
\end{align} 
So we have $\mu_{1\omega_{11}}= \dfrac{1}{2\sqrt[3]{3}}$ and from the equation \ref{p1} we have 
\begin{align}
\mu_{2\omega_{11}} =\dfrac{1}{4} \sqrt[3]{9} +\dfrac{c}{2}. 
\end{align}
We can write the energy $c$ as a function of the first root $r_1$ by 
\begin{align}
c(r_1) = \dfrac{16r_1^3 -1}{8r_1^2}. 
\end{align}
Therefore we have the torus $T_{2,1}$ for $\omega_{11}$ as
\begin{align}
T_{2,1} = \mu^{-1}(\mu_{1\omega_{11}} , \mu_{2\omega_{11}}). 
\end{align}
The critical energy value $c_{2,1}^- = -\sqrt[3]{4}$ gives us two different relation for the weight $W_2$ in the portion $\omega_{11}$. These two different relation have two different  energy cases. Namely $c\leq c_{2,1}^-$ and $c_{2,1}^- \leq c \leq -\frac{3}{2}$.

\textbf{Case 1:} Let $c \leq c_{2,1}^- =-\sqrt[3]{4}$ and the second root $r_2$ of the equation \ref{p2}. 

Since there is a tangent point of the slope $S_{2,1}$ and the graph of the equation \ref{p1}. only in $r_2$. The root $r_2$ is the length of the sides of a isosceles rightangled triangle as follows 
\begin{center} 
\includegraphics[scale=10]{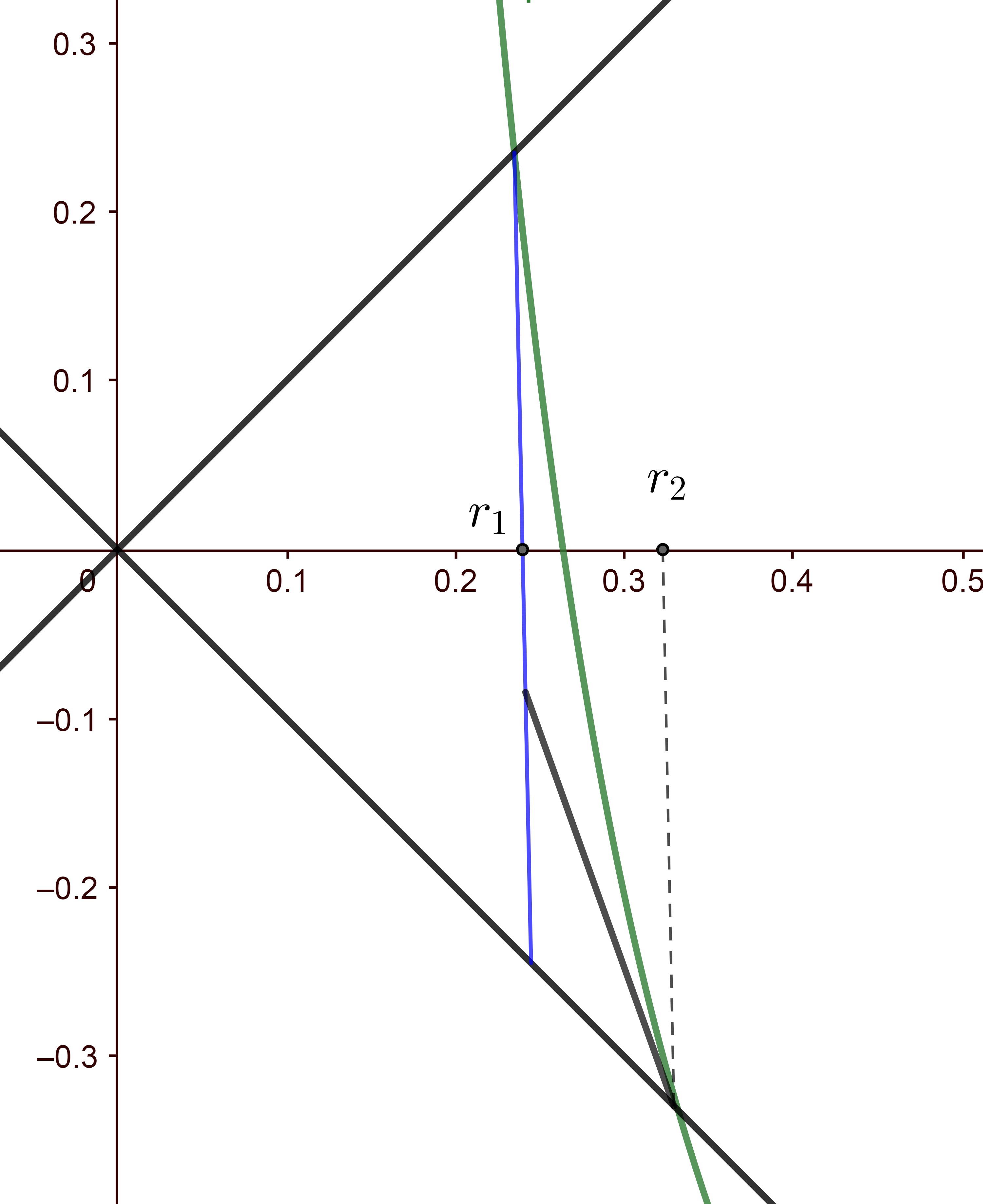}
  \captionof{figure}{$c \leq c_{2,1}^- =-\sqrt[3]{4}$} 
  \label{pic.W2-1}
\end{center}

Using the Pythagorean theorem, the diameter of this triangle is $\sqrt{2}r_2$. Therefore, considering the figure \ref{pic.W2-1}, we can write the weight $W_2$ of $\omega_{11}$ of the ECH capacities of the RKP for the case 1 as the following smooth function of $r_1$ 
\begin{align}
W_2 (r_1) & = \sqrt{2} r_2 - W_1 (r_1)= \sqrt{2} (r_2-r_1)  \\ \nonumber
& = \sqrt{2}\left[ ( - \dfrac{-1 + \sqrt{1-4^3 r_1^3}}{32r_1^2}) - r_1 \right].
\end{align}  
\textbf{Case 2:}
Let $-\sqrt[3]{4} =c_{2,1}^- \leq c \leq -\dfrac{3}{2}$. Given the figure \ref{pic.W2-2} and consider the critical point $(\mu_{1\omega_{11}} , \mu_{2\omega_{11}})$. 
\begin{center}
\includegraphics[scale=0.1]{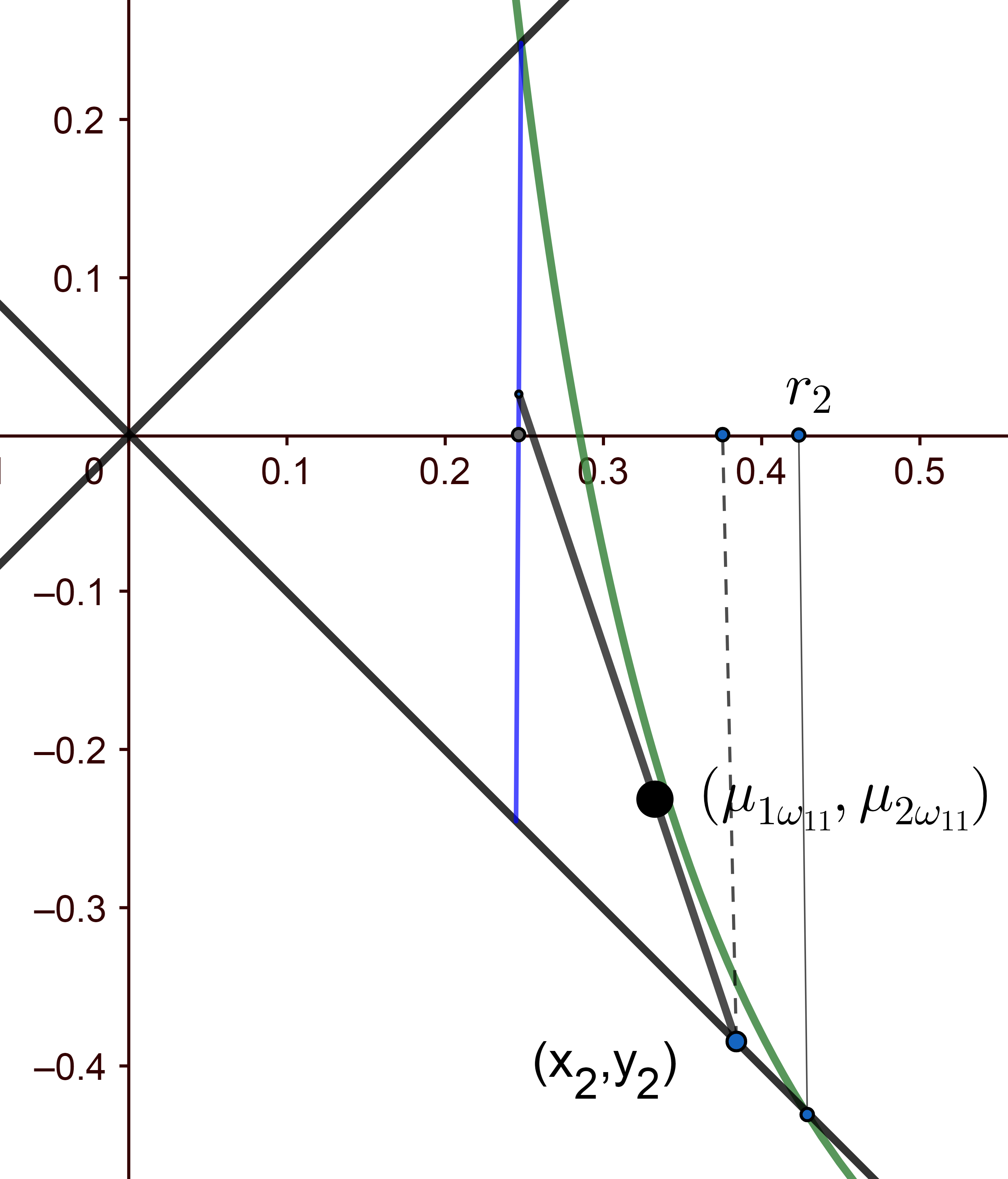}
  \captionof{figure}{$-\sqrt[3]{4} = c_{2,1}^- \leq c \leq -\frac{3}{2}$}
  \label{pic.W2-2}
\end{center}
The slope $S_{2,1}$ is tangent to the graph of the equation  \ref{p1}. in the SCTD $\mathcal{K}_c^b$ just on the critical point $(\mu_{1\omega_{11}} , \mu_{2\omega_{11}})$. Hence this tangent point is determined uniquely by the equation \ref{p1}.  

Consider the figure \ref{pic.W2-2}  and the point $(\mu_{1\omega_{11}} , \mu_{2\omega_{11}})$. We can write a line function  through  the point $(x_2,y_2)$ with the slope $-3$ as follow 
\begin{align}
y_2-\mu_{2\omega_{11}} = 3 (\mu_{1\omega_{11}} - x_2). 
\end{align}
Since $y=-x$, so the above equality becomes 
\begin{align}
-x_2 -\mu_{2\omega_{11}} = 3 (\mu_{1\omega_{11}} - x_2).
\end{align}
Therefore, $x_2$  formulated as a function of $r_1$ as 
\begin{align}
x_2(r_1) = \dfrac{3}{4 \sqrt[3]{3}} +\dfrac{\sqrt[3]{9}}{8}+\frac{1}{2} (\dfrac{16 r_1^3 -1}{16r_1^2}). 
\end{align}
Using the Pythagorean theorem, the hypotenuse of the rightangle triangle both legs have length  $x_2$ is  $\sqrt{2} x_2$. 

Now according to the figure \ref{pic.W2-2}, the weight $W_2$ of the ECH capacities of the RKP in case 2 is given as a function of  $r_1$  by 
\begin{align}
W_2(r_1) &= \sqrt{2}x_2 -W_1 (r_1) \\
&= \sqrt{2}(x_2 - r_1) \\
& = \sqrt{2} ( \dfrac{3}{4 \sqrt[3]{3}} +\dfrac{\sqrt[3]{9}}{8}+\frac{1}{2} (\dfrac{16 r_1^3 -1}{16r_1^2})-r_1). 
\end{align}
For convenience, we can write 
\begin{align}
W_2(r_1) = 
\begin{cases}
\sqrt{2}(r_2-r_1) = \sqrt{2}(-\dfrac{-1+\sqrt{1-4^3r_1^3}}{32r_1^2} - r_1),  \qquad \qquad  r_1 \leq \frac{1}{4}, \\
\sqrt{2}(\dfrac{3}{4 \sqrt[3]{3}} +\dfrac{\sqrt[3]{9}}{8}+\dfrac{1}{2} (\dfrac{16 r_1^3 -1}{16r_1^2}) -r_1), \qquad \qquad\quad \quad \qquad \frac{1}{4} \leq r_1   \leq \frac{1}{2}. 
\end{cases}
\end{align}
where $r_2=\dfrac{1}{4}$ corresponds to the energy $c=-\sqrt[3]{4}$. 

The function $W_2$ is piecewise analytic. It is also continuous at $r_1=\dfrac{1}{4}$ but is it not smooth at this point. 

\subsection{The Third Weight $W_3$: }


Here we give the computation of the third weight $W_3$ for the portion $\omega_{111}$ in the SCTD $\mathcal{K}_c^b$. 

To this deal, we need to consider the portion either $\omega_{110}$  or $\omega_{111}$, which they called even or odd respectively due to their indexes. Note that we allowed to consider the portion $\omega_{110}$  when the both cases of the second weight $W_2$ are established. But if only when  the second weight $W_2$ satisfies in the  case 2, we  can use the portion $\omega_{111}$.

Here we assume the second weight $W_2$ satisfy in the case 2 and compute the third weight $W_3$  for the portion $\omega_{111}$. 

Using the Sterb-Brocot tree, for the portion with index $111$ we can specify the node $\frac{3}{1}$ on the Hutchings CTD. On the other hand, we determine the torus $T_{3,1}$ by the critical point of the portion $\Omega_{111}$ in the CTD and  and the portion $\omega_{111}$ in SCTD $\mathcal{K}_c^b$. This torus belongs to the asteroid Hestia. Thus if we obtain the critical energy $c_{3,1}^-$ for the torus $T_{3,1}$, then we know when the asteroid Hestia appears first. 

From the relation \ref{p3} and the new tree, we have the slope $S_{3,1}=-2$ and the critical energy value 
\begin{align}
c_{3,1}^- = -\dfrac{5}{6} \sqrt[3]{6}
\end{align}
for the portion $\omega_{111}$. Using the equation \ref{p5} and \ref{p1} and the slope $S_{3,1}=-2$ gives us the critical point $(\mu_{1Hes} ,\mu_{2Hes})$ of the portion $\omega_{111}$ by 
 \begin{align}
 \dfrac{d \mu_{2}}{d \mu_1} =& - \dfrac{1}{8 \mu_{1\omega_{111}}^3 } = -2, \\
 \Longrightarrow \mu_{1\omega_{111}} =& \sqrt[3]{\frac{1}{16}}. 
 \end{align}
and 
 \begin{align}
 \mu_{2\omega_{111}} =& \dfrac{1}{16 \mu_{1\omega_{111}}^2} +\frac{1}{2}(\dfrac{16 r_1^3 -1}{8 r_1^2}) \\
 = & \dfrac{1}{16 (\sqrt[3]{\frac{1}{16}})^2} + \dfrac{16r_1^3 -1}{16r_1^2}.
 \end{align}
Using the critical point $(\mu_{1Hes} ,\mu_{2Hes})$ we can obtain the torus corresponds to the portion $\omega_{111}$ as 
\begin{align}
T_{3,1}=\mu^{-1}(\mu_{1Hes} ,\mu_{2Hes}).
\end{align}
Consider the energy values $c_{2,1}^-$ and $c_{3,1}^-$. These energies give us three different cases for the third  weight $W_3$ such that each case lives in a certain energy level. 

 \textbf{Case 1:}  Let $c \leq c_{2,1}^-$. In this case, the third weight $W_3$ of the ECH capacities of the rotating Kepler problem in the  SCTD $\mathcal{K}_c^b$ is zero. 

\textbf{Case 2:}
Let $c_{2,1}^- \leq c \leq c_{3,1}^-$. The method of this case is the same as the method of case 1 in the second weight $W_2$. Therefore, we can obtain the diameter of the isosceles rightangle triangle with the length  sides $r_2$ in the figure \ref{pic.W33-1}.
\begin{center}
\includegraphics[scale=0.2]{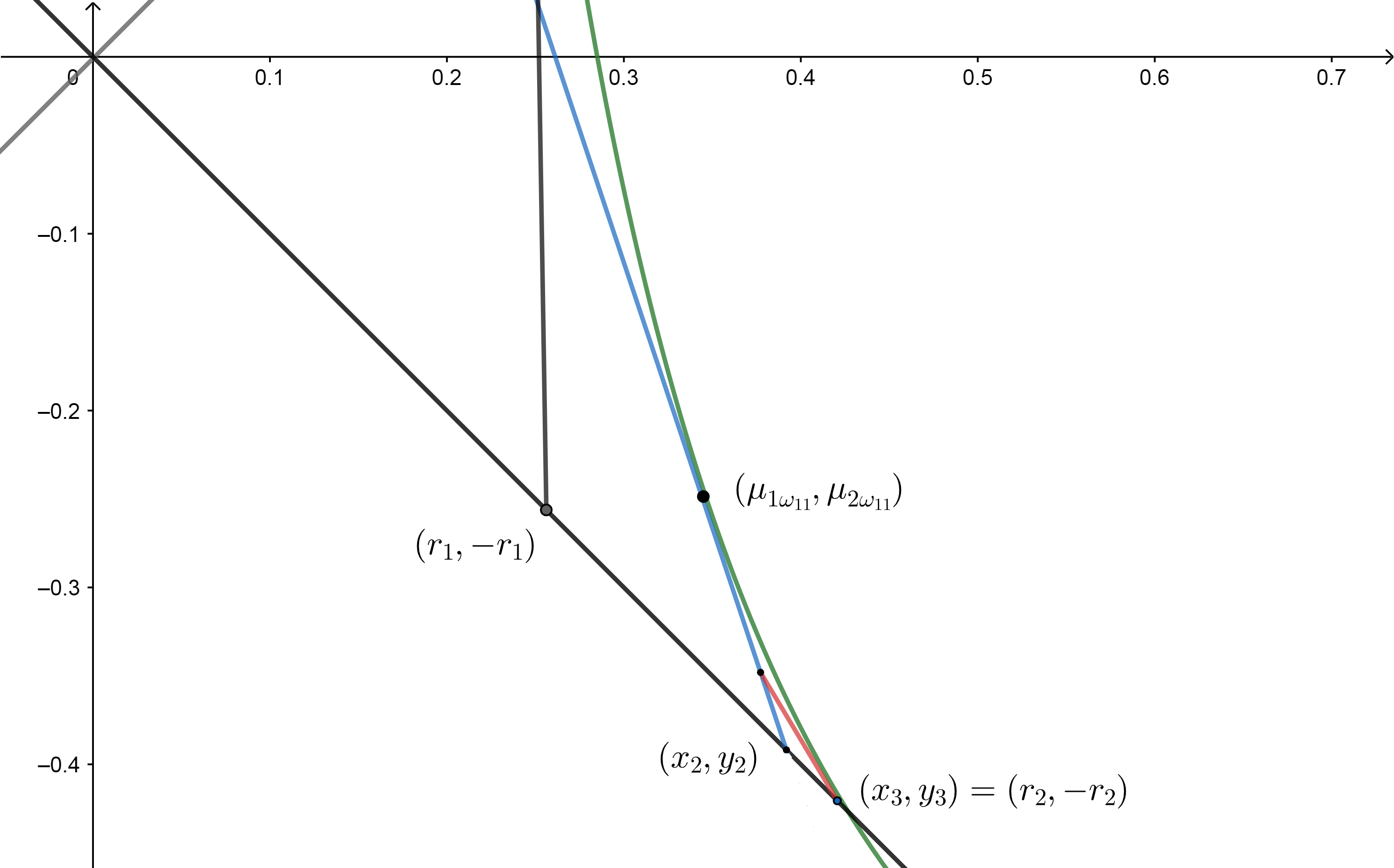}
  \captionof{figure}{Case 2 for $\omega_{111}$}
  \label{pic.W33-1}
\end{center}
Then we get the third weight $W_3$ by 
\begin{align}
W_3 = \sqrt{2}r_2 - W_2 -W_1,
\end{align}
and the weight $W_2$ is a function of $r_1$ as 
 \begin{align}
 W_3 (r_1) = \sqrt{2} ( - \dfrac{-1+\sqrt{1-4^3 r_1^3}}{32 r_1^2}) - W_2 (r_1) - \sqrt{2} r_1.
 \end{align}
\textbf{Case 3:}
Let $c_{3,1}^- \leq c \leq -\dfrac{3}{2}$. Consider figure \ref{pic.W3-2} and denote the intersection point of the slope $S_{3,1}=-2$ and line $y=-x$ by $(x_3,y_3)$ on the figure \ref{pic.W3-2}.
\begin{center}
\includegraphics[scale=15]{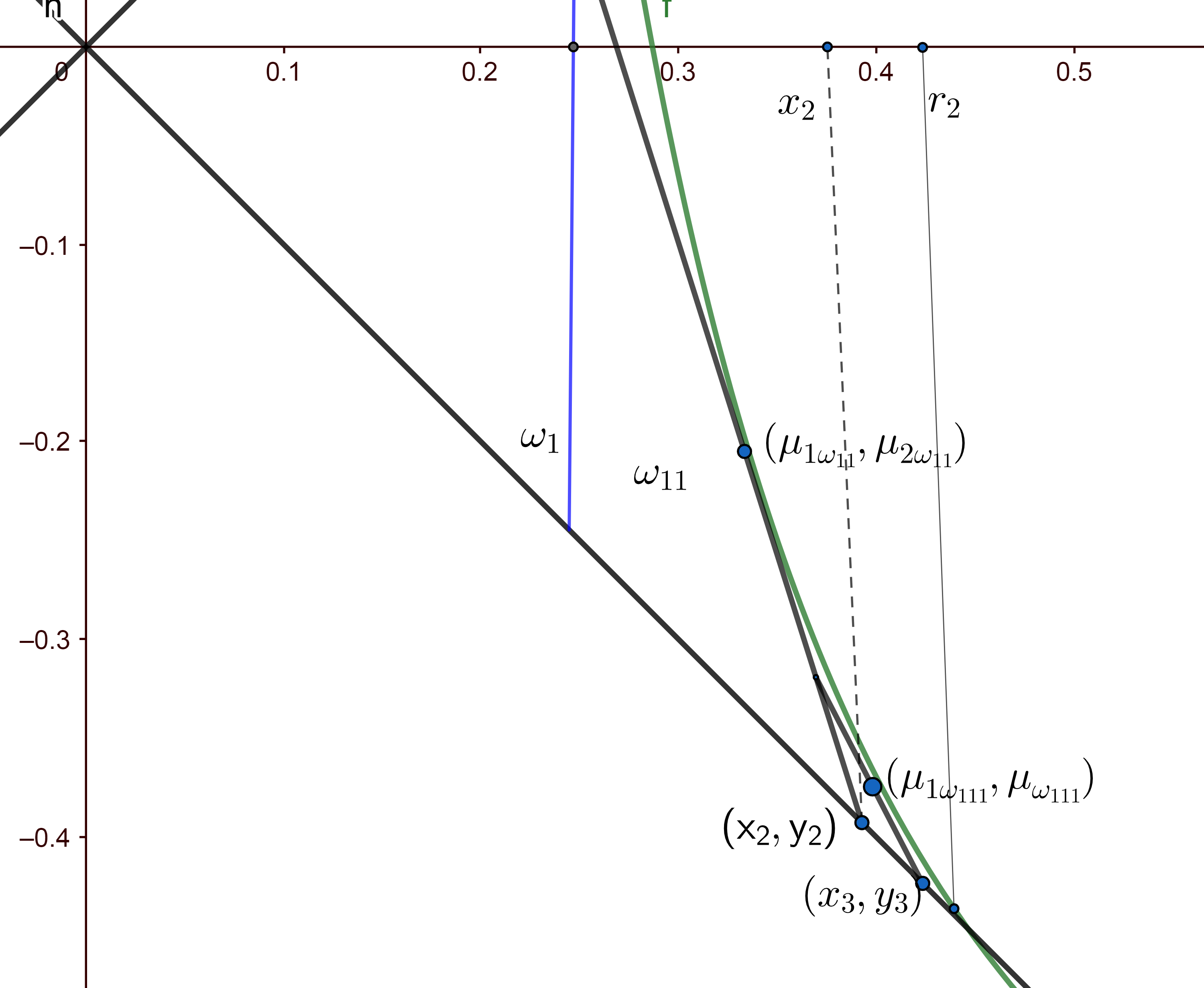}
  \captionof{figure}{$-\sqrt[3]{4} = c_{2,1}^- \leq c \leq -\frac{3}{2}$}
  \label{pic.W3-2}
\end{center}
Given the point $(\mu_{1Hes},\mu_{2Hes})$ and $(x_3,y_3)$ and write a line function of these point with the slope $S_{3,1}=-2$. Hence we have 
\begin{align}
y_3-\mu_{2Hes} = 2 (\mu_{1Hes} -x_3).
\end{align}
Since the point $(x_3,y_3)$ lives in the line $y=-x$ , so we have 
\begin{align}
x_3-\mu_{2Hes} = 2 (\mu_{1Hes} -x_3).
\end{align}
Therefore 
 \begin{align}
 x_3 = 2 \mu_{1\omega_{111}} + \mu_{2\omega_{111}} = 2 \sqrt[3]{\frac{1}{16}} +\frac{1}{16(\sqrt[3]{\dfrac{1}{16}})^2} +\dfrac{16 r_1^3 -1}{16r_1^2}.
 \end{align}
Now we follow the method of the case 2 in the second weight $W_2$  to get the third weight $W_3$ in the case 3 as a function of $r_1$ by 
 \begin{align}
 W_3(r_1) = \sqrt{2} x_3(r_1) -(W_2(r_1) +W_1(r_1)).
 \end{align} 
Finally we write the function of the third weight $W_3$ for the ECH capacities of the RKP for the energies $c \leq -\dfrac{3}{2}$ as 
 \begin{align}
 W_3(r_1) =
 \begin{cases}
 0, \qquad \qquad \qquad \qquad \qquad \qquad \qquad \qquad \qquad \qquad \qquad \qquad \qquad  r_1 \leq x_2 \\
 \sqrt{2}r_2 (r_1)-(W_2(r_1) + W_1(r_1)) = \sqrt{2} (r_2(r_1) - x_2(r_1)), \qquad  \qquad x_2 \leq r_1 \leq x_3, \\
  \sqrt{2} x_3 -(W_2(r_1) +W_1(r_1)) = \sqrt{2}(x_3(r_1) -x_2(r_1)), \qquad \quad \qquad x_3 \leq r_1 \leq r_2.
 \end{cases}
 \end{align}
\subsection{The fourth Weight $W_4$}
Here we assume the portion $\omega_{110}$ of the SCTD $\mathcal{K}_{c}^b$ and compute the fourth weight $W_4$. The portion $\omega_{110}$ is bounded by the line $x=r_1$, the graph of the equation \ref{p1} and the sloe $S_{2,1}$. We can see the portion $\omega_{110}$ in the figures \ref{pic.W4-1} and \ref{pic.W4-2}
 \begin{center}
\includegraphics[scale=14]{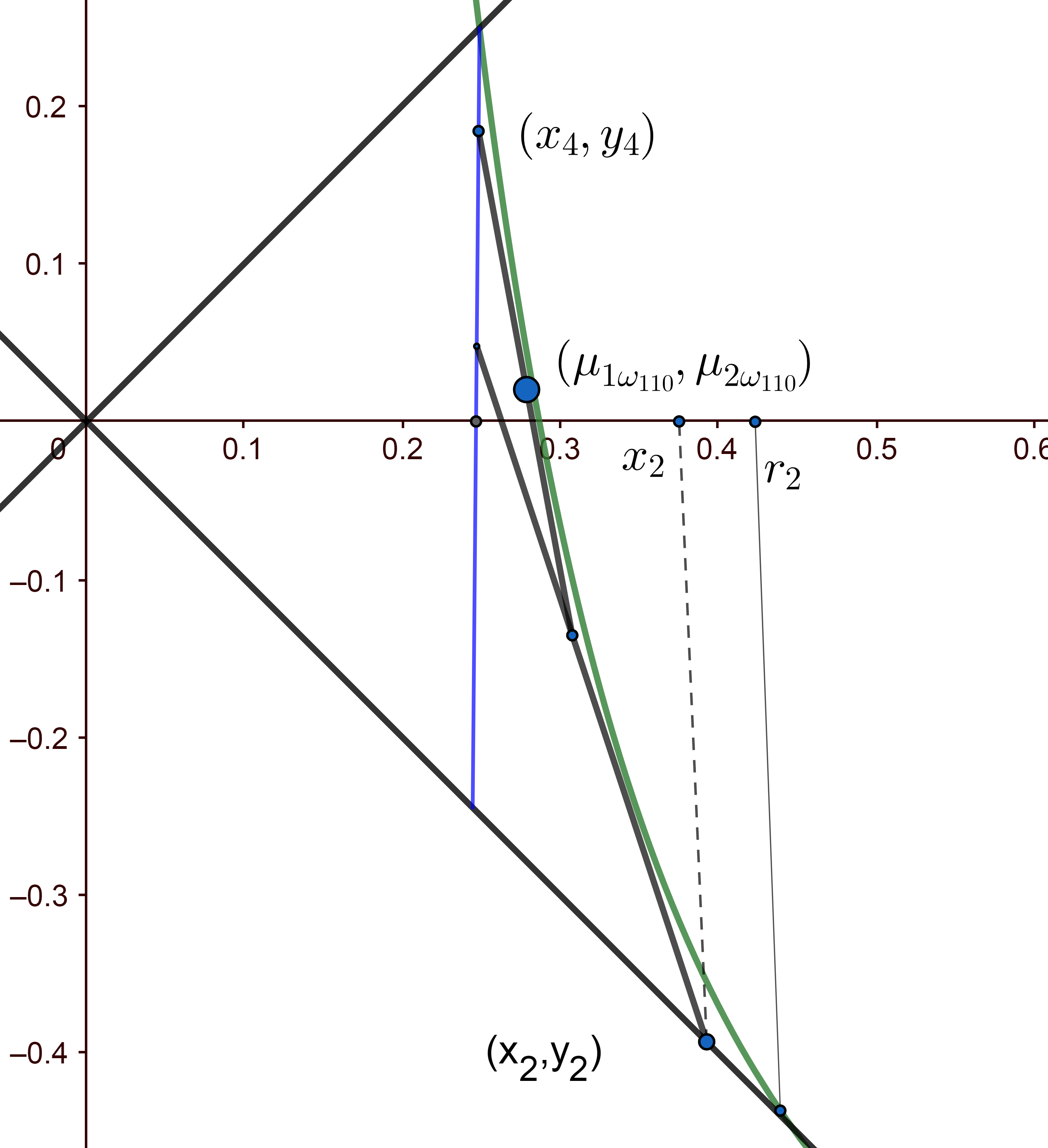}
  \captionof{figure}{$c_{3,2}^- < c_{2,1}^-$}
  \label{pic.W4-1}
\end{center}
 \begin{center}
\includegraphics[scale=0.1]{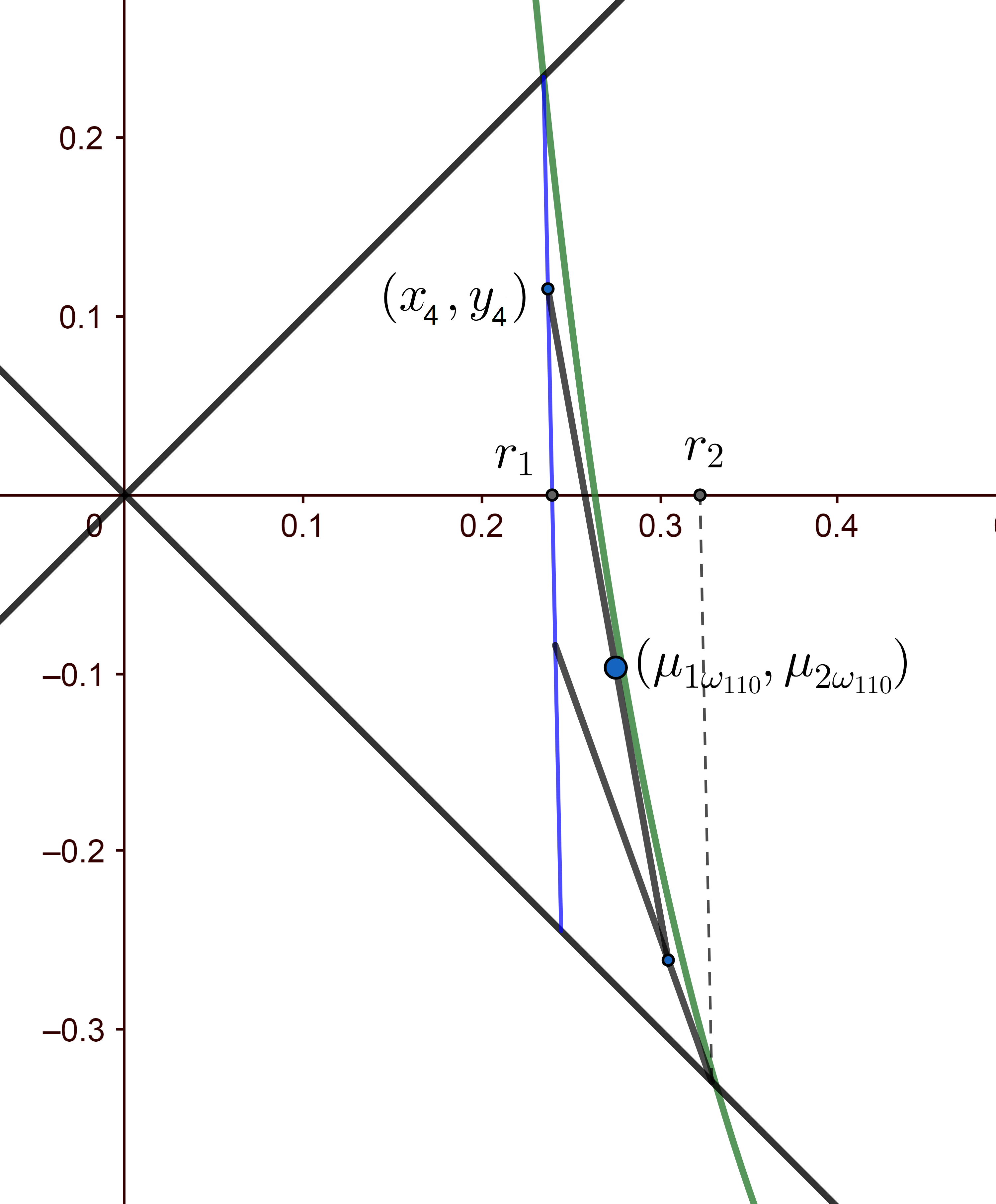}
  \captionof{figure}{$c_{3,2}^- < c_{2,1}^-$}
  \label{pic.W4-2}
\end{center}
Using the new tree gives us the slope $S_{3,2}=-5$ for the potion $\omega_{110}$ in $\mathcal{K}_c^b$. 

\begin{remark}
Unlike the third weight $W_3$ that we allowed to compute it only for the portion $\omega_{111}$ when the second weight $W_2$ satisfy in the case 2. The portion $\omega_{110}$ appears for the both cases of the second weight $W_2$. Thus we assume the portion $w_{110}$ in the SCTD $\mathcal{K}_c^b$ and obtain the weight $W_4$. 
\end{remark}

Using the Stern-Brocot tree we get the node $ \dfrac{k}{l} = \dfrac{3}{2}$ corresponds to the portion $\omega_{110}$. On the other hand, the new tree gives us the slope $S_{3,2}=-5$ for the portion $\omega_{110}$ related to the node $\dfrac{k}{l}=\dfrac{3}{2}$. 

From the following relation we get corresponding torus to this node in the Hutching CTD. 
\begin{align}
T_{3,2}=\nu^{-1} (\nu_{1\omega_{110}} , \nu_{2\omega_{110}}).
\end{align}
To compute the critical energy $c_{3,2}^-$, we have 
\begin{align}
c_{3,2}^-=-(\dfrac{1}{2} +\dfrac{3}{2})(\dfrac{3}{2})^{\frac{2}{3}} \approx -1.528768.
\end{align}
We  compute the critical point $(\mu_{1\omega_{110}} , \mu_{2\omega_{110}})$ for the fourth weight $W_4$ as 
\begin{align}
 \dfrac{d \mu_2}{d \mu_1}& = - \dfrac{1}{8 \mu_{1\omega_{110}}^3} = -5 \\
 \Longrightarrow \mu_{1\omega_{110}}&  = \sqrt[3]{\frac{1}{40}}. 
 \end{align}
and 
 \begin{align}
 \mu_{2\omega_{110}} = \dfrac{1}{16 \mu_{1\omega_{110}}^2} + \frac{1}{2} \dfrac{16 r_1^3 -1 }{ 8 r_1^2}. 
 \end{align}
Therefore, 
 \begin{align}
(\mu_{1\omega_{110}} , \mu_{2\omega_{110}}) = ( \sqrt[3]{\frac{1}{40}} , \dfrac{1}{16(\frac{1}{\sqrt[3]{40}})^2} + \dfrac{16 r_1^3-1}{16r_1^2}).
 \end{align}
From this critical point we have 
\begin{align}
T_{3,2} =\mu^{-1} (\mu_{1\omega_{110}} , \mu_{2\omega_{110}}), 
\end{align} 
in the SCTD $\mathcal{K}_c^b$. 

In the figure  \ref{pic.W4-1} and  \ref{pic.W4-2}, we named the intersection point of the slope $S_{3,2} =-5$ and the line $x=r_1$ with $(x_4,y_4)$. Since the point $(x_4,y_4)$ lives in the line $x=r_1$, we can  write the line function of the points $(r_1,y_4)$ and $(\mu_{1\omega_{110}} , \mu_{2\omega_{110}})$ with the slope $-5$ as 
\begin{align}
 y_4=& \mu_{2\omega_{110}} + 5( \mu_{1\omega_{110}}-r_1 ) \\
 y_4(r_1)=&  \mu_{2\omega_{110}} +5( \mu_{1\omega_{110}}) -5r_1  \\
 =&  \dfrac{1}{16(\frac{1}{\sqrt[3]{40}})^2} + \dfrac{16 r_1^3-1}{16r_1^2} + 5(\sqrt[3]{\dfrac{1}{40}})-5r_1 \\
 =&\dfrac{1}{16(\sqrt[3]{\dfrac{1}{40}})^2}-\dfrac{1}{16r_1^2}-4r_1 + 5(\sqrt[3]{\dfrac{1}{40}}).
 \end{align}
Finally, the fourth weight $W_4$ of the ECH capacities of the RKP is
\begin{align}
W_4(r_1) =&(r_1+y_4) - W_2(r_1) \\ \nonumber
= &( r_1 + \mu_{2\omega_{110}} +5( \mu_{1\omega_{110}}) -5r_1 ) - W_2(r_1) \\  \nonumber
=& \dfrac{1}{16(\sqrt[3]{\dfrac{1}{40}})^2}-\dfrac{1}{16r_1^2}-3r_1 + 5(\sqrt[3]{\dfrac{1}{40}}) - W_2(r_1). 
\end{align}
\begin{remark}
The necessary condition for the existence of the fourth weight $W_4$ is 
\begin{align}
s_{3,2}> \dfrac{1}{-8 r_1^3}, 
\end{align}
or equivalently 
\begin{align}
r_1< \mu_{1\omega{110}}, \qquad \qquad or \qquad \qquad r_1>\mu_{2\omega_{110}}.
\end{align}
See the figure \ref{pic.W4-3}
 \begin{center}
\includegraphics[scale=0.2]{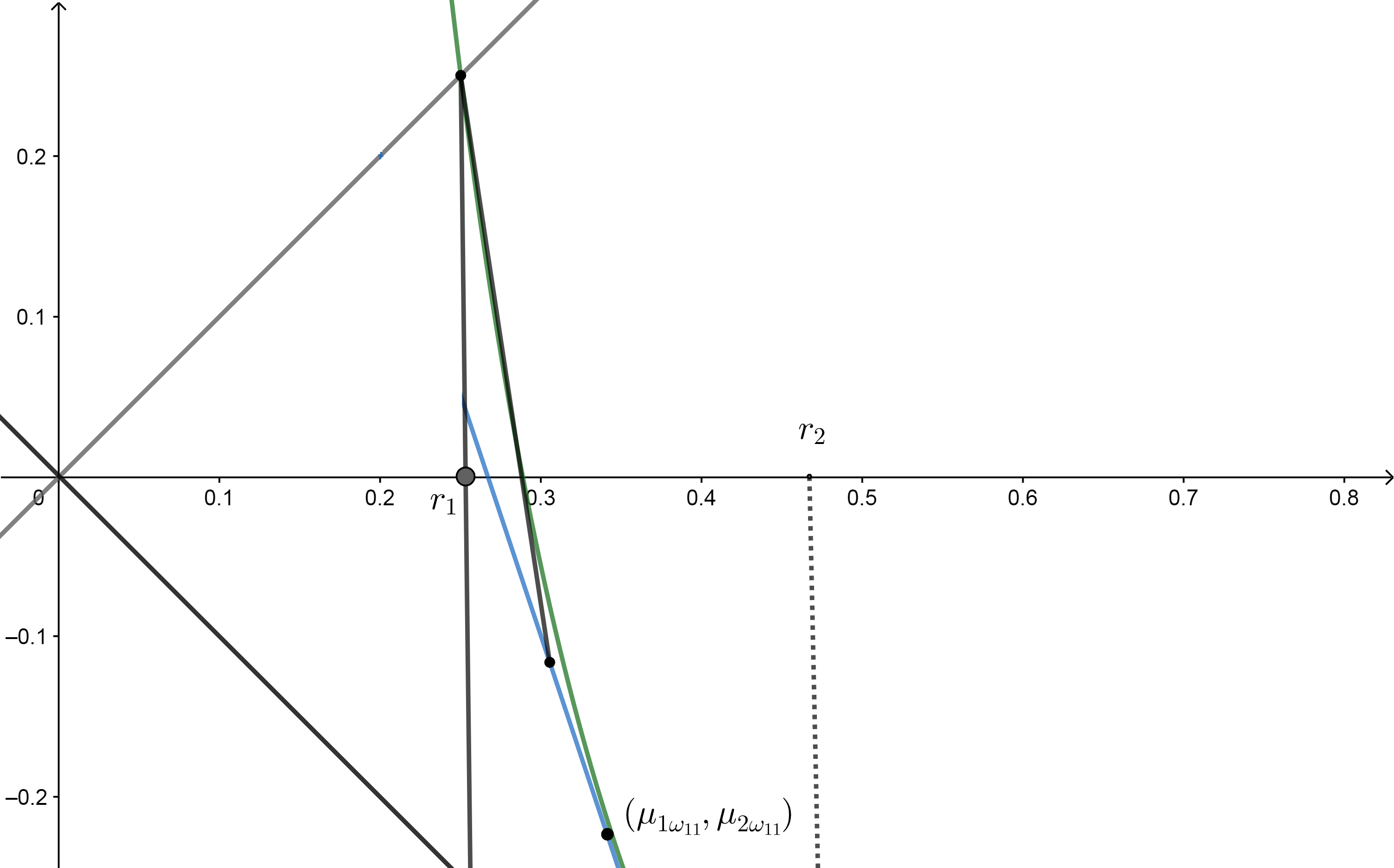}
  \captionof{figure}{$c_{3,2}^- < c_{2,1}^-$}
  \label{pic.W4-3}
\end{center}
\end{remark}
\textbf{Case 2:}
Let $c_{2,1}^- \leq  c \leq  -\dfrac{3}{2}$. Since the portion $\omega_{110}$ is not defined for $c_{2,1}^- \leq  c \leq  -\dfrac{3}{2}$. We have 
\begin{align}
W_4(r_1) = 0.
\end{align}
We abbreviate the weight $W_4$ as 
\begin{align}
W_4(r_1) =
\begin{cases}
\dfrac{1}{16(\sqrt[3]{\dfrac{1}{40}})^2}-\dfrac{1}{16r_1^2}-4r_1 + 5(\sqrt[3]{\dfrac{1}{40}}) - W_2(r_1), \qquad  c_{3,2}^- < c_{2,1}^-, \\
0,  \qquad \qquad \qquad \qquad \qquad \qquad \qquad \qquad \qquad c_{2,1}^- \leq c \leq -\frac{3}{2}.
\end{cases}
\end{align}
\begin{remark}\label{e23}
If we named  the regions $\omega_{...0}$ which is ended by zero on the SCDT $\mathcal{K}_c^b$ by the even region. Then we can generalize the above necessary condition of the existence  of the weight for  the even region as follow
\begin{align}
S_{k,l} > \dfrac{1}{-8 r_1^3},  \qquad or  \qquad
r_1 < \mu_{1\omega_{i_1 \cdots i_j}}   \qquad  or  \qquad r_1 > \mu_{2\omega_{i_1 \cdots i_j}}. 
\end{align}
where all of these three conditions are equivalent with each other, \cite{key-1}. 
\end{remark}

\subsection{The fifth weight $W_5$:}
Here we consider the portion $\omega_{1100}$ and compute the weight $W_5$ for the SCTD $\mathcal{K}_c^b$. Form the new tree we can find the slope $S_{4,3}=-7$ for this portion which  corresponds to the node $V_{1100}$ in the new tree. 

The equations \ref{p5} and \ref{p1}, we compute the critical point $(\mu_{1\omega_{1100}} , \mu_{2\omega_{1100}})$ for the portion $\omega_{1100}$ as follows
 \begin{align}
\dfrac{d \mu_2}{d \mu_1} =- \dfrac{1}{8 \mu_{1\omega_{1100}}^3} =-7  \\
\Longrightarrow \mu_{1\omega_{1100} }= \sqrt[3]{\dfrac{1}{56}}, 
\end{align}
and 
\begin{align}
\mu_{2\omega_{1100}} =& \dfrac{1}{16(\mu_{1\omega_{1100}})^2} + \dfrac{16r_1^3-1}{16r_1^2} \\ 
=& \dfrac{1}{16(\sqrt[3]{\dfrac{1}{56}})^2} - \dfrac{1}{16r_1^2}+r_1 .
\end{align}

Hence form the critical point $(\mu_{1\omega_{1100}} , \mu_{2\omega_{1100}})$ and the equation \ref{p6}
\begin{align}
T_{4,3}=\mu^{-1}(\mu_{1\omega_{1100}} , \mu_{2\omega_{1100}}). 
\end{align} 
The equation \ref{p3} gives us the critical energy value
\begin{align}
c_{4,3}^- =  (\dfrac{3}{4})(\dfrac{4}{3})^{\frac{2}{3}}.
\end{align}
This is the energy that the torus $T_{4,3}$ appears first. 

The critical energy value $c_{4,3}^-= (\frac{3}{4})(\frac{4}{3})^{\frac{2}{3}}$ gives us two different case for the weight $W_5$ such that their relations are smooth and at the critical  energy  $c_{4,3}^-= (\frac{3}{4})(\frac{4}{3})^{\frac{2}{3}}$ are continuous. 

\textbf{Case1:}
Let $c\leq c_{4,3}^-$. We follow the  method of the case 1 of the weight $W_4$ and use the slope $S_{4,3}=-7$ and the critical point $(\mu_{1\omega_{1100}} , \mu_{2\omega_{1100}})$ to obtain $y_4$ in the figure \ref{pic.W5-1}
 \begin{center}
\includegraphics[scale=24]{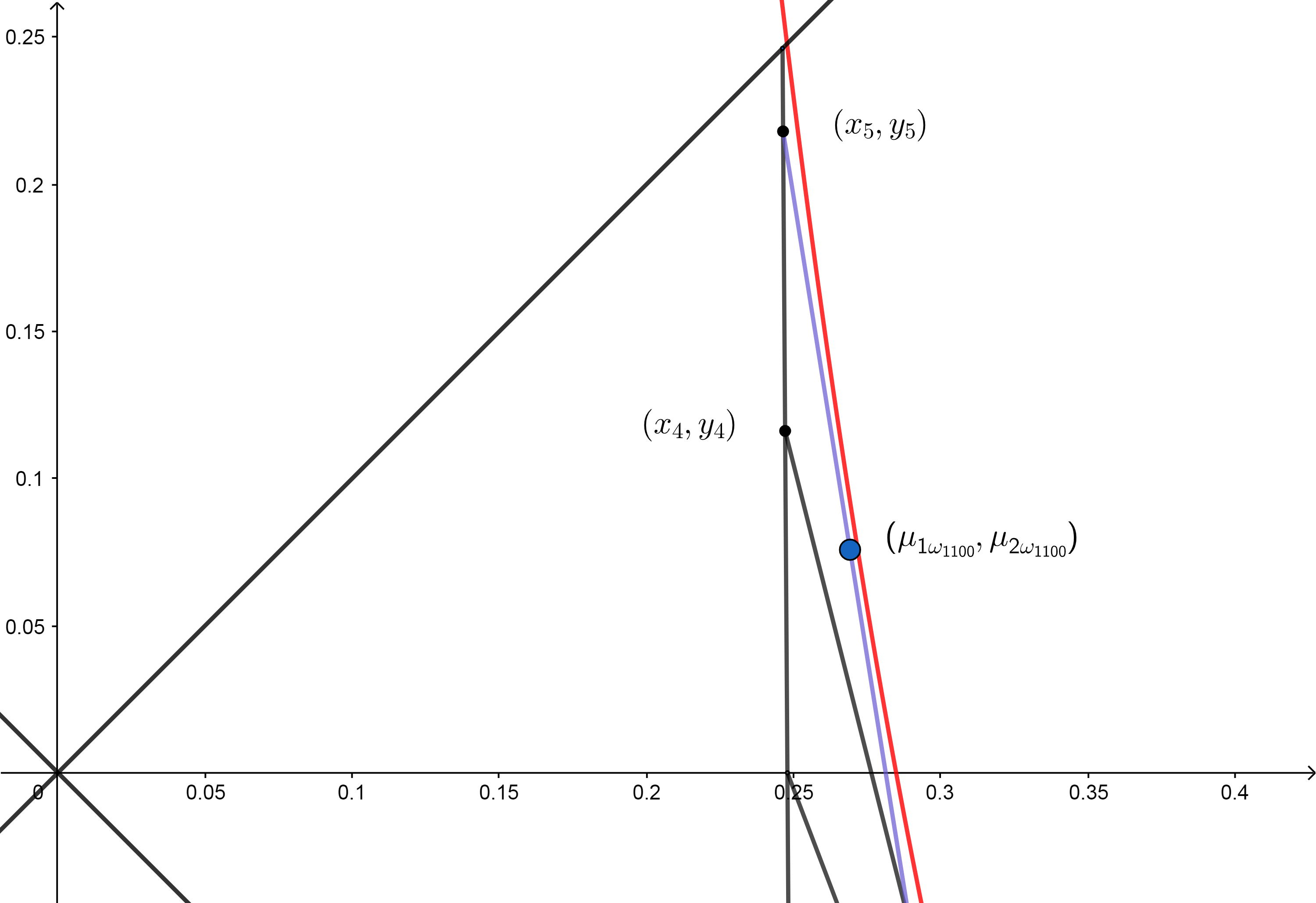}
  \captionof{figure}{The portion $\omega_{1100}$}
  \label{pic.W5-1}
\end{center}
Therefore, we can write a line function using the above notation and considering the point $(x_5,y_5)=(r_1,y_5)$ by 
\begin{align}
y_5  =&-7(r_1 - \mu_{1\omega_{1100}}) + \mu_{2\omega_{1100}} \\ \nonumber
=& -7r_1+7(\sqrt[3]{\dfrac{1}{56}}) + \dfrac{1}{16(\sqrt[3]{\dfrac{1}{56}})^2} + \dfrac{16r_1^3-1}{16r_1^2} \\\nonumber
=& -6r_1 + 7(\sqrt[3]{\dfrac{1}{56}}) + \dfrac{1}{16({\sqrt[3]{\dfrac{1}{56}}})^2} - \dfrac{1}{16r_1^2}. 
\end{align}
Using the figure \ref{pic.W5-1}, the fifth weight $W_5$ in case 1 for the SCTD $\mathcal{K}_c^b$ is 
\begin{align}
W_5(r_1) =& r_1 + y_5 -(W_2(r_1)+W_4(r_1)) \\
=& -5r_1+7(\sqrt[3]{\dfrac{1}{56}}) + \dfrac{1}{16({\sqrt[3]{\dfrac{1}{56}}})^2} - \dfrac{1}{16r_1^2} -(W_2(r_1)+W_4(r_1)). 
\end{align}
\textbf{Case 2:}
Let $c_{4,3}^- \leq c \leq -\dfrac{3}{2}$. The portion $\omega_{1100}$ is not define for the energy $c_{4,3}^- \leq c \leq -\frac{3}{2}$. Thus the fifth weight $W_5$ is equal to zero, i.e. 
\begin{align}
W_5(r_1) = 0.
\end{align}
Abbreviation of the fifth weight $W_5$ for the energy $c \leq -\dfrac{3}{2}$ is
\begin{align}
W_5(r_1)= 
\begin{cases}
 -5r_1+7(\sqrt[3]{\dfrac{1}{56}}) + \dfrac{1}{16({\sqrt[3]{\dfrac{1}{56}}})^2} - \dfrac{1}{16r_1^2} -(W_2(r_1)+W_4(r_1))\qquad \qquad c\leq c_{4,3}^-, \\
0 \qquad \qquad \qquad \qquad \qquad \qquad \qquad  c_{4,3}^- \leq c \leq -\dfrac{3}{2}.
\end{cases}
\end{align}
Note that the conditions in the remark \ref{e23} hold here. In other word, one of the following condition should hold for the portion $\omega_{1100}$, 
\begin{align}
r< \mu_{1\omega_{1100}},   \qquad  or  \qquad r_1 > \mu_{2\omega_{1100}}, \qquad or \qquad S_{4,3} > -\dfrac{1}{8r_1^3}.
\end{align}
\begin{remark}
Te weight $\omega_{i_1} ,\cdots ,\omega_{i_1 ,\cdots , i_k}$ of the SCTD $\mathcal{K}_c^b$ to for the ECH capacities of the RKP is exactly sides of isosceles rightangle triangle in the standard coordinate. 
\end{remark}
To show the above claim, we can take the weight  $\omega_{i_1, \cdots , i_k}$. There is a one-to-one corresponds between the Stern-Brocot tree and the new tree and also the correspondence condition holds for the the portions $\Omega_{i_1,\cdots ,i_k}$ and $\omega_{i_1,\cdots ,i_k}$via the new tree  in the SCTD $\mathcal{K}_c^b$. That means, the node $V_{i_1 ,\cdots , i_k}$ in the Stern-Brcot tree gives us the portion $\Omega_{i_1,\cdots ,i_k}$ in the Hutchings CTD is equivalent to the portion $\omega_{i_1,\cdots ,i_k}$ in the SCTD $\mathcal{K}_c^b$. 

If we rotate the portion $\omega_{i_1,\cdots ,i_k}$ by 45 degree in counter clockwise direction then multiplying it to
$
\left[ \begin{array}{ccc}
1 & 0 \\
1 & 1 
\end{array} \right]
 \in SL_2(\mathbb{Z})$
 and  
  or 
 $
\left[ \begin{array}{ccc}
1 & 1 \\
0 & 1 
\end{array} \right]
 \in SL_2(\mathbb{Z})$,
  we will obtain an isosceles rightangle triangle with slope $-1$  corresponding the the portion $\Omega_{i_1,\cdots ,i_k}$  in the Hutchings CTD.  Note that the slope of the portion has a one-to-one relation with the node $V_{i_1,\cdots ,i_k}$ in the Stern-Brocot tree. Hence sometimes we need to do the above multiplication   several times to get the isosceles right-angles triangle with the slope -1. 
\section{The Integrals of the Region}
From the definition of the ECH capacities we need to give an order to the weights from the biggest one to lowest one. To obtain this order, we compute the area of the regions that the weights are defined on those areas. 

We compute the area of the $\omega_1$ for the first weight $W_1$. This area is an isosceles rightangle triangle. Thus we have 
\begin{align}\label{e12}
Area( \omega_1) = \frac{1}{2} (\sqrt{2} r_1)^2 = r_1^2. 
\end{align}
The area of the rest part of the SCTD $\mathcal{K}_c^b$ is computing by the following integral,
\begin{align}\label{e11}
Area(\mathcal{K}_c^b - \omega_1) & = \int_{r_1}^{r_2} \mu_2 - (-\mu_1) d \mu_1 \\ \nonumber
&=\int_{r_1}^{r_2} \dfrac{1}{16\mu_1^2} + \frac{c}{2}-(-\mu_1) d\mu_1 \\ \nonumber
&= -\dfrac{1}{16\mu_1}+ \frac{1}{2}c\mu_1 + \frac{1}{2}\mu_1^2 |_{r_1}^{r_2}\\ \nonumber
&=\left(-\dfrac{1}{16(-\dfrac{-1 + \sqrt{1-4^3r_1^3}}{32r_1^2})} + \dfrac{1}{2}(\dfrac{16r_1^3 -1}{8r_1^2})(-\dfrac{-1+\sqrt{1-4^3r_1^3}}{32r_1^2}) + \dfrac{1}{2}(-\dfrac{-1+\sqrt{1-4^3r_1^3}}{32r_1^2})^2 \right)  \\ \nonumber
& \qquad - \left(-\dfrac{1}{16r_1}+\dfrac{1}{2}(\dfrac{16r_1^3 -1}{8r_1^2})r_1 +\dfrac{1}{2}r_1^2 \right) \\ \nonumber 
&= \dfrac{2r_1^2}{-1+\sqrt{1-4^3r_1^3}} +\dfrac{(1-\sqrt{1-4^3r_1^3})-2(-2+24r_1^3)}{32r_1}  \\ \nonumber
& \qquad \qquad + \dfrac{4(-1+\sqrt{1-4^3r_1^3}) + 2-2\sqrt{1-4^3r_1^3} -4^3r_1^3}{2 \times 32^2 r_1^4} \\ \nonumber
&= \dfrac{2 r_1^2}{-1+\sqrt{1-4^3r_1^3}} + \dfrac{5-\sqrt{1-4^3r_1^3} -48r_1^3}{32r_1} + \dfrac{-2+2\sqrt{1-4^3r_1^3} -64r_1^3}{2 \times 32^2 r_1^4}. \nonumber
\end{align}
\begin{example}\label{e9}
If we let $c =-\frac{3}{2}$. Then we have $r_1 = \frac{1}{4}$ and $r_2 = \frac{1}{2}$ and also 
\begin{align}
Area(\omega_1) = & (\frac{1}{4})^2 = \frac{1}{16} \\
Area(\mathcal{K}_c^b - \omega_1) =& \int_{r_1=\frac{1}{4}}^{r_2=\frac{1}{2}} \mu_2 + \mu_1 d \mu_1 = \frac{1}{32}. 
\end{align}
\end{example}

\begin{theorem}\label{e21}
Let $c \leq -\dfrac{3}{2}$. We have the following inequality, 
\begin{align}
\mathcal{F}(r_1)  = \dfrac{Area(\mathcal{K}_c^b - \omega_1)} {Area(\omega_1)} \leq \dfrac{1}{2} \qquad\qquad\qquad \forall r_1 \in [0,\dfrac{1}{4}].
\end{align}
\end{theorem}
\begin{proof} 
From the relations \ref{e11} and \ref{e12}, we have the following identities 
\begin{align}
Area(\mathcal{K}_c^b - \omega_1) = & \dfrac{2 r_1^2}{-1+\sqrt{1-4^3r_1^3}} + \dfrac{5-\sqrt{1-4^3r_1^3} -48r_1^3}{32r_1} + \dfrac{-2+2\sqrt{1-4^3r_1^3} -64r_1^3}{2 \times 32^2 r_1^4}, \\ \nonumber
Area(\omega_1) = & r_1^2. 
\end{align}
We take these two identities and compute the following relation
\begin{align}
\dfrac{Area(\mathcal{K}_c^b - \omega_1)}{Area(\omega_1)} &=\dfrac{\dfrac{2 r_1^2}{-1+\sqrt{1-4^3r_1^3}} + \dfrac{5-\sqrt{1-4^3r_1^3} -48r_1^3}{32r_1} + \dfrac{-2+2\sqrt{1-4^3r_1^3} -64r_1^3}{2 \times 32^2 r_1^4}}{r_1^2} \\ \nonumber 
&= \dfrac{2}{-1+\sqrt{1-4^3r_1^3}} + \dfrac{5-\sqrt{1-4^3r_1^3} - 48r_1^3}{32r_1^3} +\dfrac{-1+\sqrt{1-4^3r_1^3} -32r_1^3}{32^2r_1^6} \\ \nonumber 
& = \dfrac{2}{-1+\sqrt{1-4^3r_1^3}} + \dfrac{5-\sqrt{1-4^3r_1^3}}{32r_1^3} - \dfrac{3}{2} +\dfrac{-1 +\sqrt{1-4^3r_1^3}}{32^2r_1^6 }- \dfrac{1}{32r_1^3 }. \nonumber 
\end{align}
For simplicity, we define $R:= r_1^3$. Thus we have 
\begin{align}\label{e13}
\dfrac{Area(\mathcal{K}_c^b - \omega_1)}{Area(\omega_1)} = & \dfrac{2}{-1+\sqrt{1-4^3R}} + \dfrac{5-\sqrt{1-4^3R}}{32R} - \dfrac{3}{2} +\dfrac{-1 +\sqrt{1-4^3R}}{32^2R^2 }- \dfrac{1}{32R } \\ \nonumber 
=& \dfrac{2}{-1+\sqrt{1-4^3R}} + \dfrac{4-\sqrt{1-4^3R}}{32R} - \dfrac{3}{2} +\dfrac{-1 +\sqrt{1-4^3R}}{32^2R^2 } \\ \nonumber 
=& \dfrac{2}{-1+ \sqrt{1-4^3R}} +\dfrac{-1-\sqrt{1-4^3R}}{-1-\sqrt{1-4^3R}} + \dfrac{4- \sqrt{1-4^3R}}{32R} + \dfrac{-1+\sqrt{1-4^3R}}{32^2R^2} -\dfrac{3}{2} \\ \nonumber
=& \dfrac{3-2\sqrt{1-4^3R}}{32R} + \dfrac{-1+\sqrt{1-4^3R}}{32^2R^2} -\dfrac{3}{2}. 
\end{align}
Now we take the first derivative of the equation \ref{e13} respect to $R$. Thus we have 
\begin{align}
\dfrac{d(\dfrac{Area(\mathcal{K}_c^b - \omega_1)}{Area(\omega_1)})}{dR} &= 
\dfrac{-64R-3\sqrt{1-4^3R} +2}{32 R^2\sqrt{1-4^3R} } + \dfrac{48R + \sqrt{1-4^3R} - 1}{512R^3\sqrt{1-4^3R} } \\ \nonumber 
&= \dfrac{-1024 R^2 - 48R\sqrt{1-4^3R} + 80R + \sqrt{1-4^3R} -1}{512R^3 \sqrt{1-4^3R}}.
\end{align}
If we put the nominator of the above equation equal to zero, then we can get the zeros of the nominator at the points $R=0$ and $R=\dfrac{1}{64} = \dfrac{1}{4^3}$. From the definition, we have $r _1 = \sqrt[3]{R} = \dfrac{1}{4}$. Therefore the function $\mathcal{F}$ on its domain $(0 , \dfrac{1}{4}]$ is monotone increasing. 
Note that the function $\mathcal{F} $ take its maximal value at the point $\dfrac{1}{4}$, i.e. $\mathcal{F}(\dfrac{1}{4}) = \dfrac{1}{2}$, which we have already obtained in the example \ref{e21}. 
\end{proof}
\begin{example}
In this example we are going to compute the ECH capacities for the RKP by using the SCTD $\mathcal{K}_c^b$ for some $k$ when the energy $c=-\dfrac{3}{2}$. 

First if we use the equations of weights which are introduced  in this chapter, we can get the following values
\begin{align}
W_1(c) = &   \sqrt{2}r_1 \approx  0.353554 \\ \nonumber 
W_2(r_1) \approx  &      0.219247 \\ \nonumber 
W_3(r_1) \approx  &   0.0502325 \\ \nonumber 
W_4(r_1) \approx  &   0.223766 \\ \nonumber 
W_5(r_1) \approx  &    0.0514663 \\ \nonumber 
\end{align}
Note that for the energy $c=-\dfrac{3}{2}$, we have $r_1= \dfrac{1}{4}$, $r_2=\dfrac{1}{2}$  and form examples \ref{e9},  we know 
\begin{align}
Area(\mathcal{K}_c^b) &= \dfrac{3}{32} \\ \nonumber 
Area(\omega_1) &= \dfrac{1}{16} \\ \nonumber 
Area(\mathcal{K}_c^b - \omega_1)& = \dfrac{1}{32} \\ \nonumber 
Area(\omega) & \approx 0.01501571682 .\\ \nonumber 
\end{align}
And also Theorem \ref{e21} says that $W_1(c)$ is the first weight of the ECH capacities of the RKP. 

Therefore the above computations give us the following order of  the weights $W_1, \cdots , W_5$ and $W_j$, 
\begin{align}
W_1 > W_4>W_2>W_5>W_3 \geq W_j \qquad  \forall j \in \mathbb{N}  \text{ and } j \geq 6.
\end{align}
Now consider the inequality 
\begin{align}
d^2 +d \leq 2k
\end{align}
Than we can have the following table. 
\begin{center}
\captionof{table}{ECH capacities for $c=-\dfrac{3}{2}$}
\label{Table}
\begin{tabular}{|c|c|c|}
\hline
 Rank & The ECH cap. for $\mathcal{K}_c^b$ & The ECH cap. for $c=-\dfrac{3}{2}$ \\
\hline
$c_1(\mathcal{K}_c^b)$ & $W_1 $ & $ 0.353554 $ \\ \hline
$c_2(\mathcal{K}_c^b)$ & $W_1+W_4 = c_1+W_4 $ & $0.57732$ \\ \hline
$c_3(\mathcal{K}_c^b)$ & $2W_1=2c_1 $ & $0.707108$ \\ \hline
$c_4(\mathcal{K}_c^b)$ & $2W_1 +W_4=c_3+W_2$ & $0.930874$ \\ \hline
$c_5(\mathcal{K}_c^b)$ & $2W_1+W_4+W_2=c_4+W_2 $ & $1.150121$ \\ \hline
$c_6(\mathcal{K}_c^b)$ & $2W_1+2W_4 =2c_2$ & $1.15464$ \\ \hline
$c_7(\mathcal{K}_c^b)$ & $3W_1+W_4=3c_1+W_4 $ & $1.284428$ \\ \hline
$c_8(\mathcal{K}_c^b)$ & $3W_1+W_4+W_2=c_7+W_2 $ & $1.503675$ \\ \hline
$c_9(\mathcal{K}_c^b)$ & $3W_1+2W_4=c_7+W_4 $ & $1.508194$ \\ \hline
$c_{10}(\mathcal{K}_c^b)$ & $3W_1+2W_4+W_2 =c_9+W_2$ & $1.727441$ \\ \hline
$c_{20}(\mathcal{K}_c^b)$ & $5W_1+W_4+W_2+W_5 $ & $2.2622493$ \\ \hline
\end{tabular} 

\end{center}
\end{example}

\textit{Center for Mathematical Analysis, Geometry and Dynamical Systems,
Instituto Superior T\'ecnico, Universidade de Lisboa,
Av. Rovisco Pais, 1049-001 Lisboa, Portugal.}

\textit {\author{ E-mail address:  \href{mailto:aminmohebbi@tecnico.ulisboa.pt}{aminmohebbi@tecnico.ulisboa.pt} }}


\end{document}